\documentclass[a4paper]{article}

\usepackage{natbib}
\usepackage[dvips]{graphicx}
\usepackage{amssymb,amsmath,xcolor}
\usepackage{amsxtra, amsthm,authblk}
\usepackage{upref}
\usepackage{geometry}
\usepackage{tikz}
\usetikzlibrary{arrows.meta}

\usepackage[colorlinks,
            linkcolor=blue,
            anchorcolor=green,
            citecolor=blue
            ]{hyperref}

\numberwithin{equation}{section}

\DeclareMathOperator{\dist}{dist}

\newcommand{\E}{\mathbb{E}}
\newcommand{\A}{\mathbb{A}}
\newcommand{\R}{\mathbb{R}}
\newcommand{\C}{\mathbb{C}}
\newcommand{\N}{\mathbb{N}}
\newcommand{\Z}{\mathbb{Z}}
\renewcommand{\P}{\mathbb{P}}
\newcommand{\W}{\mathbb{W}}
\newcommand{\cA}{\mathcal{A}}
\newcommand{\cB}{\mathcal{B}}
\newcommand{\cE}{\mathcal{E}}
\newcommand{\cF}{\mathcal{F}}
\newcommand{\cG}{\mathcal{G}}
\newcommand{\cH}{\mathcal{H}}
\newcommand{\cK}{\mathcal{K}}
\newcommand{\cL}{\mathcal{L}}
\newcommand{\cV}{\mathcal{V}}
\newcommand{\harm}{\mathrm{harm}}
\newcommand{\diam}{\mathrm{diam}}
\newcommand{\vol}{\mathrm{vol}}
\newcommand{\Var}{\mathrm{Var}}
\newcommand{\Cov}{\mathrm{Cov}}
\newcommand{\Area}{\mathrm{Area}}

\newtheorem{theorem}{Theorem}[section]
\newtheorem{lemma}[theorem]{Lemma}
\newtheorem{proposition}[theorem]{Proposition}
\newtheorem{definition}[theorem]{Definition}
\newtheorem{problem}[theorem]{Problem}
\newtheorem{assumption}[theorem]{Assumption}

\theoremstyle{remark}
\newtheorem{claim}[theorem]{Claim}

\title{Internal diameter bounds for DLA on random planar maps under invariance principles}
\author{Deven Mithal\thanks{University of Chicago}}
\date{}
\begin{document}

\maketitle

\begin{abstract}
   We establish sublinear bounds for the internal diameter of External Diffusion Limited Aggregation (DLA) clusters on random planar maps. We develop two conditional DLA growth results for deterministic planar maps: the weaker condition yielding sublinear growth for DLA, and the stronger yielding polynomially sublinear growth for DLA. We then verify the strong assumption a.s.\ for the $\gamma$-mated CRT maps, and the weak condition a.s.\ for ergodic scale-free random environments, establishing the corresponding a.s.\ bounds.
\end{abstract}

\tableofcontents

\section{Introduction}\label{sec:intro}
External Diffusion Limited Aggregation (DLA) is a discrete-time random growth process on a graph, consisting of a sequence of increasing subgraphs, in which new edges and vertices are added at each step according to the harmonic measure from a target vertex or $\infty$. While a number of variants of DLA have been studied, rigorous results regarding the process have been scarce. This paper is principally concerned with establishing bounds on the size and shape of DLA on random planar geometries as time tends to $\infty$. Two notions of size can be considered: 1) external diameter, in which the DLA cluster is viewed as a subset of the ambient graph, and diameter of the cluster is measured with respect to the ambient graph metric, and 2) internal diameter, in which the DLA cluster is viewed as a graph in its own right, and the diameter of that graph is measured. The latter notion is more informative of the shape of the DLA cluster. Indeed, consider a spiral curve on a graph versus a radially expanding ball. These increasing subgraphs may have the same external diameter, but from the perspective of the internal diameter these are polar opposites.

While there does exist limited literature on the external diameter of DLA on random planar maps, to our knowledge no sublinear upper bounds have been established for the internal diameter. The purpose of this paper is to derive these bounds on the internal diameter of DLA. The primary tool used in this paper is \textit{quenched} control of the random walk on the random planar map, i.e.\ control of the random walk \textit{conditional} on the random environment. Such properties have been established for large classes of random planar maps (see e.g. \cite{gwynne2021invarianceprincipleergodicscalefree}), and so the results of this paper are expected to generalize to a number of models.

As the inputs to our arguments are almost sure quenched statements, we typically work with a deterministic underlying geometry, and indeed our main results (Theorems~\ref{thm:sublinear} and \ref{thm:poly-sublinear}) are stated as black boxes for deterministic planar maps. Transfer to random planar map models then follows by a simple application of the deterministic theorem to the relevant random planar map model on the appropriate event (Section~\ref{subsec:random-transfer}).

We deal with two notions of quenched control of random walk: 1) almost sure quenched convergence of random walk to Brownian motion (Assumption~\ref{assumption:rw-conv}), and 2) almost sure quenched mesoscopic control of random walk (Assumption~\ref{assumption:mesoscopic}). The first assumption is sufficient to establish that the internal diameter of the DLA cluster is $o(n)$, and the second (stronger) assumption allows one to prove that the internal diameter of the DLA cluster is $O(n^{1 - \zeta})$ with $\zeta > 0$. For our principal example, which is the $\gamma$-mated continuum random tree (CRT) map (Definition~\ref{def:mated-CRT}), the exponent $\zeta$ will be related to the Liouville quantum gravity (LQG) parameter $\gamma$.

We now define DLA, and follow the presentation in \cite{gwynne2021external}. Below, $\cV(G)$ denotes the vertices of the graph $G$; see Section~\ref{sec:notation} for a more detailed discussion of notation.

\begin{definition}[External DLA]\label{def:dla}
    Let $G$ be an infinite, connected graph, and let $v_0, v_* \in \cV(G)$. External DLA from $v_0$ to $v_*$ is a sequence $(K_n)_{n \in \N_0}$ of random subgraphs of $G$ with $K_0 = \{v_0 \}$, and with $K_n$ defined as follows:
    \begin{enumerate}
        \item If $v_* \in K_{n - 1}$, terminate the process, i.e. $K_{n + k - 1} = K_{n - 1}$ for all $k \in \N$.
        \item Otherwise, run a simple random walk from $v_*$, conditioned to hit $K_{n - 1}$, and consider the edge $e_n$ first traversed from $G \setminus K_{n - 1}$ into $K_{n - 1}$. Define $K_n$ as the union of $K_{n - 1}$, $e_n$ and the endpoint $v_n$ of $e_n$ in $G \setminus K_{n - 1}$.
    \end{enumerate}
    External DLA targeted at $\infty$ is the same process with the harmonic measure from $v_* \in G$ replaced with harmonic measure from $\infty$, where such a measure is well defined.
\end{definition}

Unless otherwise stated, we work with recurrent rooted graphs $(G,o)$, DLA starts from $v_0 = o$ and is targeted at $\infty$, and harmonic measure from $\infty$ is assumed to be unique and well-defined. A rooted set is a vertex subset containing $o$. For a fixed deterministic graph, almost sure statements about DLA refer to its growth randomness.

There are several notions of diameter for which DLA has been studied. \cite{gwynne2021external} concerns the diameter of $K_n$ considered as a subset of $G$, which is the so-called \textit{external diameter}, and is denoted $\diam(K_n; G)$. This paper concerns the \textit{internal diameter}, which considers $K_n$ as a graph in its own right, and is denoted $\diam(K_n; K_n)$. As this paper only concerns the internal diameter, we write $\diam(K_n)$ for short.

\begin{definition}[DLA growth exponent]\label{def:growth-exponent}
    The (internal) growth exponent of DLA from $o$ to $\infty$ on a connected rooted graph $(G,o)$ is defined as the limit
    \begin{equation*}
        \lim_{n \to \infty} \frac{\log \diam(K_n)}{\log n}
    \end{equation*}
    when it exists.
\end{definition}

We now proceed to state our main results, and the properties on which they rely.

\subsection{Property 1: Quenched convergence of random walk to Brownian motion}\label{subsec:prop-1}

First, we require a few definitions, beginning with an embedding of a graph into $\C$. For the definition of a proper embedding in the plane, we follow the presentation of \cite[Section 2.4]{hutchcroft2019uniform}.

\begin{definition}[Proper embedding in the plane]\label{def:embedding}
    Define a proper embedding $\Phi$ of a graph $G$ into the plane as a locally finite topological embedding of vertices $\cV(G)$ and edges $\cE(G)$ into $\C$.
\end{definition}

Subsequently, we refer to the object in Definition~\ref{def:embedding} as an \textit{embedding}. We normalize the embedding so that $\Phi(o) = 0$. A random walk on $G$ is mapped into $\C$ under the embedding in the obvious way; we refer to this as the \textit{embedded walk}. For finite curves parametrized by $[0,1]$, set
\begin{equation*}
  d_F(\gamma_1, \gamma_2)
  := \inf_{\phi} \sup_{t \in [0,1]}
    |\gamma_1(t) - \gamma_2(\phi(t))|
\end{equation*}
where $\phi$ ranges over increasing homeomorphisms of $[0,1]$. The function $d_F$ induces a metric on curves viewed modulo parametrization, with curves of zero $d_F$-distance considered equivalent \cite{AB99}.

For $R,L > 0$, let $Y^{x;R,L}$ be the embedded walk started from $x$ and stopped at its first exit from $B_{LR}(0)$ (which may be in the interior of an edge). Subsequently, we will choose $L$ as a radius multiplier and $R$ as the non-normalized scale for the embedded graph. A normalized scale $L$ corresponds to the radius $r = L R$ before rescaling. Let $B^{z;L}$ be standard planar Brownian motion started from $z$ and stopped on leaving $B_L(0)$. We write $\cL$ for law and $\cL_x$ for law with starting point $x$. We can now state the first of the aforementioned invariance principles in the form of Assumption~\ref{assumption:rw-conv}.

\begin{assumption}[Random walk convergence hypothesis]\label{assumption:rw-conv}
    For every $0 < C < L < \infty$, we have
    \begin{equation*}
      \sup_{x:\,|\Phi(x)| \leq CR} d_{\textrm{Pr}}(\cL(R^{-1}Y^{x;R,L}), \cL(B^{\Phi(x)/R;L})) \to 0 \text{ as } R \to \infty
    \end{equation*}
    where $d_{\textrm{Pr}}$ is the Prokhorov distance induced by $1 \wedge d_{F}$.
\end{assumption}

The above assumption states that the rescaled embedded random walk converges in law to Brownian motion (modulo parametrization) uniformly over starting points in compact subsets of $\C$. Using Assumption~\ref{assumption:rw-conv} we can establish an analogue of Kesten's Beurling estimate on $\Z^d$. Below, $\harm_A^\infty$ is the harmonic measure on the edges out of $A$ viewed from $\infty$, and $\cE_A(v)$ is the set of edges out of $A$ incident to $v \in A$; see Section~\ref{sec:notation} for a more precise discussion.

\begin{theorem}[Weak Beurling estimate]\label{thm:beurling-weak}
    Under Assumption~\ref{assumption:rw-conv}, we have
    \begin{equation*}
      \lim_{m \to \infty} \sup_{\substack{A \ni o \text{ finite, connected}\\ |A| \geq m}} \, \sup_{v \in A}
      \harm_A^\infty(\cE_A(v)) = 0.
    \end{equation*}
\end{theorem}

Using Theorem~\ref{thm:beurling-weak}, we can prove sublinearity of the internal diameter of the DLA cluster. We state this result as a theorem rather than a corollary, as DLA is the focus of this paper.

\begin{theorem}[Sublinear growth]\label{thm:sublinear}
    Let $(K_n)$ be DLA started from the root $o$ and targeted at $\infty$. Under Assumption~\ref{assumption:rw-conv}, we have
    \begin{equation*}
      \diam(K_n) = o(n)
    \end{equation*}
    almost surely.
\end{theorem}

The conclusion also holds for DLA started from any other vertex $v$: adjoining a fixed path from $o$ to $v$ transfers the harmonic measure bound from sets $A$ with $o \in A$ to sets $v \in A$, up to a finite multiplicative constant depending on the path (see Equation~\eqref{eq:root-comparison}).

\subsection{Property 2: Quenched mesoscopic control of random walk}\label{subsec:prop-2}
The assumption and results of this subsection require slightly more information than that of Assumption~\ref{assumption:rw-conv}. Our planar map, rather than being presented by merely an embedding, is given by a \textit{cellular embedding}. We combine the presentation from \cite[Definition 1.6]{gwynne2021tutteembeddingmatedcrtmap} and \cite{gwynne2019harmonicfunctionsmatedcrtmaps}.

\begin{definition}[Cellular embedding in the plane]\label{def:cellular-embedding}
    A cellular embedding into $\C$ consists of:
    \begin{enumerate}
        \item A locally finite collection of compact connected subsets of $\C$ with nonempty interiors, denoted by $\cH$, such that $\cup_{A \in \cH} A = \C$ and $\mathrm{Leb}(A \cap B) = 0$ for all distinct $A, B \in \cH$.
        \item A planar graph $G$, consisting of vertices $\cV(G)$ and unweighted edges $\cE(G)$.
        \item A bijection $\phi$ from $\cV(G)$ to $\cH$ such that $a \sim b$ only if $\phi(a) \cap \phi(b) \neq \emptyset$.
    \end{enumerate}
    We write $H_x := \phi(x)$ and normalize so that the root cell $H_o$ contains $0$.
\end{definition}

For $D \subset \C$ and $A \subset \cV(G)$, write $\cV_{\cH}(D) := \{x \in \cV(G): H_x \cap D \neq \emptyset \}$ and $U_A^{\cH} := \bigcup_{x \in A}H_x$. A path is said to  disconnect two vertex sets if its vertices intersect every path in the graph joining them.

We choose to present the mesoscopic result in the context of a cellular embedding because it is slightly more natural for the $\gamma$-mated CRT map application, but the subsequent assumption and theorems can be readily adapted to the planar embedding data of Assumption~\ref{assumption:rw-conv}.

The above definition is implicit in \cite{gwynne2019harmonicfunctionsmatedcrtmaps}, and is a special case of \cite[Definition 1.15]{gwynne2021invarianceprincipleergodicscalefree}. The more general definition, which is called a \textit{cell configuration}, has in place of items 2 and 3 in Definition~\ref{def:cellular-embedding} a symmetric adjacency relation $\sim$ with $A \sim B$ only if $A \neq B$ and $A \cap B \neq \emptyset$, and a symmetric function $\mathfrak{c}: \{(A, B): \, A, B \in \cH \text{ and } A \sim B \} \to (0, \infty)$ representing conductance. The graph structure can be encoded by $\mathfrak{c}(A, B) = |\{e \in \cE(G): \, e \text{ is between } \phi^{-1}(A) \text{ and } \phi^{-1}(B) \}|$. We leave the adjacency and conductances implicit in the notation $\cH$, but will discuss the more general notion of cell configuration in Section~\ref{sec:ergodic}, in which we sketch a proof for an analogue of Theorems~\ref{thm:beurling-weak} and \ref{thm:sublinear} with more general edge weights (called conductances).

We can now introduce the mesoscopic invariance principle assumption. This assumption involves three functions which depend on the non-normalized scale $R$:
\begin{enumerate}
    \item $m(R)$, which is the mesoscopic scale lower bound at non-normalized scale $R$.
    \item $e(R)$, which is the net approximation error at non-normalized scale $R$.
    \item $d(R)$, which is the maximal cell diameter at non-normalized scale $R$.
\end{enumerate}

\begin{assumption}[Mesoscopic control hypothesis]\label{assumption:mesoscopic}
    There are constants $p \in (0,1)$ and $\sigma > 0$, a minimal non-normalized scale $R_0 < \infty$, and functions
    $$
      m, e, d: [R_0, \infty) \to (0, \infty)
    $$
    such that the following hold.
    \begin{enumerate}
    \item $m(R) \leq R^{1 - \sigma}$ and $\frac{e(R) + d(R)}{m(R)} \to 0$.

    \item Every cell intersecting $B_{4R}(0)$ has Euclidean diameter at most $d(R)$.
    
    \item For every $R \geq R_0$ there exist finite $e(R)$-nets $\mathcal{Z}_R \subset B_{2R}(0)$ and $\mathcal{Q}_R \subset [m(R), R/100]$ such that for every $(z,r) \in \mathcal{Z}_R \times \mathcal{Q}_R$ and every $x \in \cV_{\cH}(\partial B_r(z))$ a simple random walk started from $x$ has probability at least $p$ to disconnect $\cV_{\cH}(\partial B_{3r/4}(z))$ from $\cV_{\cH}(\partial B_{5r/4}(z))$ before exiting $\cV_{\cH}(\A_{r/2, 3r/2}(z))$, with the union of cells visited up to and including exit contained in $\A_{r/3, 5r/3}(z)$.

    \item $D(A) := \diam(U_A^{\cH}) \geq c_{\vol}|A|^{\nu}$ for some $c_{\vol}, \nu > 0$, for every sufficiently large finite connected $A$ containing $o$.
    \end{enumerate}
\end{assumption}
We now briefly describe the above assumption in words. Item 1 controls the size of the mesoscopic scale lower bound relative to the non-normalized scale, net error and diameter, item 2 states the cell diameter control enforced by $d(R)$, item 3 provides a lower bound on the disconnection probability for each fixed scale and center in the nets, and 4 provides control of the diameter of the cellular embedding of a vertex subset in terms of its cardinality which enables the passage to intrinsic graph geometry from embedded graph geometry.

The use of nets in item 3 is to enable a union bound argument which provides a quantitative upper bound on the failure probability of the above assumption for random planar maps, which is necessary to apply Borel-Cantelli and conclude mesoscopic control at all sufficiently large non-normalized scales. In principle, one could strengthen item 3 to hold for all $(z, r) \in B_{2R}(0) \times [m(R), R/100]$, which would allow one to weak the limit condition in item 1 to $d(R)/m(R) \to 0$, but this becomes more difficult to verify for e.g. the $\gamma$-mated CRT maps.

As with Assumption~\ref{assumption:rw-conv}, we can prove a Beurling-type estimate using Assumption~\ref{assumption:mesoscopic}, albeit a stronger one than that of Theorem~\ref{thm:beurling-weak}.

\begin{theorem}[Strong Beurling estimate]\label{thm:beurling-strong}
    Under Assumption~\ref{assumption:mesoscopic}, there exists $\theta > 0$ and a finite threshold $N = N(G,o)$ such that
    \begin{equation*}
      \harm_A^\infty(\cE_A(v)) \leq |A|^{-\theta}
    \end{equation*}
    for all finite rooted connected $A$ with $|A| > N$ and all $v \in A$.
\end{theorem}

Using Theorem~\ref{thm:beurling-strong}, we establish the polynomially sublinear bounds on the internal diameter of DLA, again stated as a theorem.

\begin{theorem}[Polynomially sublinear growth]\label{thm:poly-sublinear}
    Let $(K_n)$ be DLA started from the root $o$ and targeted at $\infty$. Under Assumption~\ref{assumption:mesoscopic}, there exists a deterministic $\zeta > 0$ such that
    \begin{equation*}
      \diam(K_n) = O(n^{1 - \zeta})
    \end{equation*}
    almost surely.
\end{theorem}

As with Theorem~\ref{thm:sublinear}, this result can be transferred to an arbitrary initial vertex. The stronger mesoscopic control allows one to upper bound the growth exponent by $1 - \zeta$.

\subsection{Transfer to random planar maps}\label{subsec:random-transfer}

The primary random planar map model of interest in this paper will be the so-called mated CRT map, which we now define. We follow the presentation of \cite[Section 1.1]{gwynne2019harmonicfunctionsmatedcrtmaps}.

\begin{definition}[Mated CRT map]\label{def:mated-CRT}
    Fix $\gamma \in (0, 2)$ and let $(L, R)$ be a two sided two-dimensional Brownian motion such that
    \begin{equation*}
        \Var(L_t) = \Var(R_t) = |t| \quad \text{ and } \quad \Cov(L_t, R_t) = -|t| \cos \bigg(\frac{\pi \gamma^2}{4} \bigg)
    \end{equation*}
    and normalized with $L_0 = R_0 = 0$. The $\epsilon$-mated CRT map for $(L, R)$ is the random graph with vertices $\epsilon \Z$, and where vertices $x_1, x_2 \in \epsilon \Z$ with $x_1 < x_2$ connected by an edge if and only if
    \begin{equation*}
        \bigg(\inf_{t \in [x_1 - \epsilon, x_1]}L_t \bigg) \vee \bigg(\inf_{t \in [x_2 - \epsilon, x_2]}L_t \bigg) \leq \inf_{t \in [x_1, x_2 - \epsilon]} L_t
    \end{equation*}
    or the same condition holds for $R$. If both conditions hold for $x_1 + \epsilon < x_2$, then the vertices are connected by two edges. We denote this random planar map, rooted at $o = 0$, by $\cG^{\epsilon}$. We observe that by identifying $x \in \epsilon \Z$ with $x/\epsilon \in \Z$ the law of $\cG^{\epsilon}$ does not depend on $\epsilon$ by Brownian scaling.
\end{definition}

When the parameter $\gamma$ is obvious or unimportant, we refer to the above graph simply as the mated CRT map. Note that the above definition does not immediately specify a planar embedding of $\cG^{\epsilon}$. However, due to the relationship between the mated CRT map, Schramm-Loewner evolution (SLE), and LQG, it does come equipped with a canonical cellular embedding; this is presented in Section~\ref{sec:crt-lqg}.

By an application of Theorem~\ref{thm:poly-sublinear} to the mated CRT map a.s., one can achieve an upper bound on the DLA growth exponent (Definition~\ref{def:growth-exponent}). As random planar maps are the principal concern of this paper, rather than deterministic ones, we record the consequences of Theorems~\ref{thm:sublinear} and \ref{thm:poly-sublinear} in the following theorem.

\begin{theorem}[DLA on the $\gamma$-mated CRT]\label{thm:dla-mated-crt}
    Fix $\gamma \in (0, 2)$. Then there exists a $\zeta = \zeta(\gamma) > 0$ such that DLA started from the root $o = 0$ and targeted at $\infty$ on the $\gamma$-mated CRT map satisfies $\diam(K_n) = O(n^{1 - \zeta})$ a.s.
\end{theorem}

Using \cite{gwynne2021invarianceprincipleergodicscalefree} one can make a rather strong and general claim: growth of DLA on an arbitrary ergodic scale free random environment (in particular, those obeying the hypothesis of \cite[Theorem 1.5]{gwynne2021invarianceprincipleergodicscalefree}) is a.s.\ sublinear. We now define the general class of ergodic scale-free random environments, which we call \textit{good embedded lattices}, to which our theory applies, following the presentation of \cite[Section 1.1]{gwynne2021invarianceprincipleergodicscalefree}.

\begin{definition}\label{def:ergodic}
Let $\mathcal{M}$ be a random, infinite, locally finite planar, undirected graph (with loops and multiple edges allowed) with conductances $\mathfrak c: \cE(\mathcal{M}) \to (0, \infty)$, embedded so that each edge is a simple curve with zero Lebesgue measure, edges intersect only at their endpoints, each compact subset of $\C$ intersects at most finitely many edges, and each connected component of $\C \setminus \mathcal{M}$ has compact closure. We call this class of graphs embedded lattices, and denote it by $\mathbb{M}$. A good embedded lattice satisfies:
\begin{enumerate}
\item (Translation invariance modulo scaling) For every nonnegative measurable $F: \mathbb{M} \times \C^2 \to [0, \infty)$ on the space of doubly-marked embedded lattices such that
$$
  F(a(\mathcal{M} - z), a(u - z), a(v - z)) = a^{-2}F(\mathcal{M}, u, v), \qquad a > 0, \quad z,u,v \in \C
$$
we have
$$
  \E \bigg[\int_\C F(\mathcal{M}, 0, w)\, dw \bigg] = \E \bigg[\int_\C F(\mathcal{M}, w, 0)\, dw \bigg].
$$

\item (Ergodicity modulo scaling) Every measurable event invariant under all transformations of the form $\mathcal{M} \mapsto a(\mathcal{M} - z)$ for $a > 0$ and $z \in \C$ has probability $0$ or $1$.

\item (Finite expectation) Faces are the closures of the connected components of $\C \setminus \mathcal{M}$. For $z \in \C$, write $H_z$ for a face of $\mathcal{M}$ containing $z$, chosen arbitrarily if there is more than one, and denote the collection of faces by $\cF \mathcal{M}$. For distinct faces, write $H \sim H'$ if they share an edge, and $\mathfrak c(H,H')$ for the sum of the conductances of their shared edges. Additionally, define
\begin{equation*}
    \pi(H) := \sum_{\substack{H' \in \cF \mathcal{M} \\ H \sim H'}} \mathfrak c(H,H') \quad \text{ and } \quad \pi^*(H) := \sum_{\substack{H' \in \cF \mathcal{M} \\ H \sim H'}} \mathfrak c(H,H')^{-1}
\end{equation*}
For the origin-containing face $H_0$ we have
$$
  \E \bigg[\frac{\diam(H_0)^2}{\Area(H_0)}(\pi(H_0) + \pi^*(H_0)) \bigg] < \infty.
$$
\end{enumerate}
\end{definition}

The main invariance principle pertaining to a good embedded lattice, stated as Proposition~\ref{prop:ergodic} in the preliminaries, is the convergence of the random walk on the embedded lattice to a \textit{non-isotropic} Brownian motion \cite{gwynne2021invarianceprincipleergodicscalefree}. Our methods for isotropic Brownian motion adapt to this setting in a relatively straightforward manner to produce sublinear internal diameter bounds for DLA.

Here, DLA is defined by replacing the simple random walk with the conductance weighted random walk (see Definition~\ref{def:conductance-weighted-walk}), assuming harmonic measure from infinity is unique.

\begin{theorem}[DLA on an ergodic scale-free random environment]\label{thm:dla-ergodic-env}
    Let $\mathcal{M}^*$ be the dual of a good embedded lattice $\mathcal{M}$ rooted at $o = H_0$, and let $(K_n)$ be External DLA on $\mathcal{M}^*$ started from $o$ and targeted at $\infty$ on $\mathcal{M}^*$. Then $$\diam(K_n) = o(n)$$ a.s.
\end{theorem}

The almost sure conclusions of Theorems~\ref{thm:dla-mated-crt} and \ref{thm:dla-ergodic-env} are joint in the environment and DLA growth. A proof sketch of this result is provided in Section~\ref{sec:ergodic}. Note that in both Theorems~\ref{thm:dla-mated-crt} and \ref{thm:dla-ergodic-env}, the required control of the random walk is almost sure convergence of the law of the random walk conditioned on the environment. For a number of random planar maps (RPMs), such as spanning-tree-weighted maps \cite{holden2026tutte}, only convergence in probability of the conditionals is available, and one may wonder if the techniques of this paper might be extended to this weakened assumption. In Section~\ref{sec:open-probs}, we present an example which possesses this weaker property, but for which DLA in fact has a.s.\ linear internal diameter growth. This demonstrates the tightness of the arguments of this paper.

\subsection{Related work}
The literature on random growth processes is broad and diverse, and in particular a number of distinct setups and questions for DLA have been addressed, as well as for continuum analogues and candidate scaling limits. We now briefly review this work, and group these contributions by whether they deal with a discrete or continuum growth process, and a random or deterministic underlying geometry.

\paragraph{Discrete growth on deterministic geometry.} Witten and Sander introduced DLA \cite{witten1981diffusion, witten1983diffusion} to model processes of accumulation in nature. Prior rigorous work especially relevant to this paper is the work of Kesten on hitting probabilities and so-called long arms \cite{kesten1987hitting,kesten1987long,kesten1990upper}, which produces Beurling-type bounds on the harmonic measure of the boundary of a connected lattice subset, uniformly over the particular boundary vertex and subsets of a minimal size, and uses these bounds to control the likelihood that the next accumulated particle is adjacent to the most recently accumulated one. Subsequent work has extended the technique of Kesten to various other graph geometries \cite{benjamini2008diffusion, BenjaminiYadinSpectralGap, benjamini2017upper, lawler2004beurling}. This line of work heavily exploits the use of the harmonic measure in DLA.

DLA belongs to a more general class of models called the Dielectric Breakdown Models (DBM) \cite{niemeyer1984dbm, hastings2001fractal, MJ02, MJB08, hastings1998laplacian, carleson2001aggregation, RZ05, NorrisTurnerHL}, in which the harmonic measure is replaced by an $\eta$'th power of the harmonic measure; DBM includes models such as the Eden model, which is of independent interest \cite{eden1961two, hammersley1965first, Ric73, cox1981some, bubeck2015asymptotic, curien2017scaling, curien2019first}. Recent work achieves $\eta$-dependent arm bounds for DBM on $\mathbb Z^2$ and $\mathbb Z^3$, recovering the classical DLA bounds at $\eta = 1$, by departing to an extent from the technique of Kesten \cite{losev2025long}. Other variants of DLA include one-dimensional long-range DLA \cite{amir2016one}, stationary half-plane DLA \cite{ProcacciaZhangStationary,ProcacciaYeZhangStationary}, and Internal DLA \cite{lawler1992internal,jerison2012logarithmic,asselah2013sublogarithmic}. The latter model is closely related to the Abelian sandpile \cite{bou2021convergence, bou2022dynamic, bou2024shape, bou2024integer, levine2010scaling}. Both the work on Beurling-type bounds and on long-arms control the internal diameter of DLA.

\paragraph{Continuum growth on deterministic geometry.} A number of continuum aggregation growth models on deterministic geometries have been studied, which provide candidate scaling limits for DLA on deterministic geometries; these may be regarded as ``semi-discrete'' models, as they are discrete-time models but accumulate continuum segments. The Hastings-Levitov (HL) model is a family of conformal aggregation models \cite{hastings1998laplacian, NorrisTurnerHL}, which for a particular choice of parameter is expected to be related to DLA \cite[Section 1.2.3]{ms-qle}; some scaling limit results have been established for these models, as well as for the broader class of growth models called Aggregate Loewner Evolution (ALE) \cite{NorrisSilvestriTurnerALE, SolaTurnerViklund, HiggsALE, berger2022growth}, which encompasses the HL model. However, rigorous relationships between HL/ALE and DLA have yet to be achieved.

\paragraph{Discrete growth on random geometry.} An early contribution to the study of DLA on random geometries is \cite{benjamini2017upper}. For the uniform infinite spanning-tree-weighted random planar map (UITM), Gwynne and Pfeffer prove that the DLA cluster has external diameter $n^{2/d + o(1)}$, for an exponent $d$ related to volume of the UITM \cite{gwynne2021external}; their work exploits heavily the relationship between DLA and the Loop-Erased Random Walk (LERW) on the UITM. On the mated CRT maps, which are studied in this paper, shape theorems for internal DLA and the divisible sandpile have been established~\cite{bou2024internal}.

\paragraph{Continuum growth on random geometry.} Miller and Sheffield introduced Quantum Loewner Evolution (QLE) as a candidate scaling limit for the DBM models on an RPM converging to $\gamma$-LQG \cite{ms-qle}. Their construction is motivated by the relationship between the $\eta$-DBM model and the $\eta$-Laplacian Random Walk (LRW), which generalizes that between DLA and LERW on the UITM. It is expected that DLA on the UITM converges jointly to QLE$(2,1)$ on a $\sqrt{2}$-LQG surface. The so-called ``quantum natural time'' $\operatorname{QLE}(8/3,0)$ was constructed and used to construct the $\sqrt{8/3}$-LQG metric \cite{ms-equivalence,ms-equiavlence-2}. Ang and the author subsequently construct quantum natural time QLE for a range of parameters including the $(8/3, 0)$ construction of \cite{ms-equivalence}, establish the qualitative phases of the model, as well as stationarity and Markov properties of the unexplored surface. Their time parametrization is conjectured to arise from the time parametrization of the discrete DBM models, but thus far, aside from the Eden model on uniform RPMs (see \cite[]{ang2026quantum}), convergence of DBM models on RPMs to QLE on LQG is not established.

\subsection{Outline}
The proof strategy of both Theorems~\ref{thm:sublinear} and \ref{thm:poly-sublinear} is identical: bound the maximal harmonic measure (from $\infty$) of $\cE_{K_n}(v)$ over any $v \in K_n$. Indeed, this is identified as the core mechanism for internal diameter bounds by Gwynne and Pfeffer in \cite[Section 8]{gwynne2021external}. To do this, one needs to identify the size of annulus $\mathbb{A}$ required such that the probability that the random walk, upon entering $\mathbb{A}$, disconnects the relevant vertex sets before exiting $\mathbb{A}$, possesses a positive lower bound. The larger the number of disjoint annuli which can be identified, the stronger the upper bound on harmonic measure, which is what produces the differing conclusions of Theorems~\ref{thm:sublinear} and \ref{thm:poly-sublinear}. 

In Section~\ref{sec:prelims}, we present preliminaries on the mated CRT map required to verify Assumption~\ref{assumption:mesoscopic}. In Section~\ref{sec:arg-template} we prove a bound on the internal diameter of DLA conditional on an annular intersection probability via a harmonic measure bound. In Sections~\ref{sec:pf-sublinear} and \ref{sec:pf-poly-sublinear} we use the blackbox theorem of Section~\ref{sec:arg-template} to establish Theorems~\ref{thm:sublinear} and \ref{thm:poly-sublinear} under assumptions~\ref{assumption:rw-conv} and ~\ref{assumption:mesoscopic} respectively. In Section~\ref{sec:mated-verification}, we verify Assumption~\ref{assumption:mesoscopic}. In Section~\ref{sec:ergodic} we sketch the proof of Theorem~\ref{thm:dla-ergodic-env}. Lastly, in Section~\ref{sec:open-probs} we discuss open problems related to DLA on RPMs, including the mated CRT and UITM, and present a counterexample illustrating an obstruction to a straight-forward application of known, weaker invariance principles to the study of DLA.

\subsection{Notation and terminology}\label{sec:notation}
Write $\N = \{1,2, \ldots \}$ and $\N_0 := \N \cup \{0\}$. For $K \in \N$, write $[K] := \{1, \ldots, K\}$. We identify a subgraph with its vertex set when using vertex-set notation, and $|A|$ denotes vertex cardinality. Graph distance is denoted by $\dist(\cdot, \cdot;G)$; for subsets of $\C$, $\dist$ and $\diam$ are Euclidean. Edges are unoriented, and parallel edges are distinct. Simple random walk chooses uniformly among incident edge ends, counting each self-loop twice. We denote by $\A_{r,s}(z)$ the Euclidean annulus $B_s(z) \setminus B_r(z)$ for $r < s$ and $z \in \C$. When $z = 0$, we drop the center and write $\A_{r,s}$. For $v \in G$, define the neighborhood $N(v) := \{x \in \cV(G): \, x \sim v \} \cup \{ v \}$, and for a subset $A \subset \cV(G)$ with $v \in A$, define the external neighborhood $S_A(v) := \{x \in \cV(G) \cap A^c: \, x \sim v \}$. Additionally, we let $\cE_A(v) := \{e \in \cE(G): \, e \text{ joins }v \text{ to }S_A(v) \}$, i.e. edges with one end on $v$ and the other outside of $A$. For finite nonempty $A$ and $x \notin A$, $\harm_A^x$ is the law of the edge traversed when a random walk from $x$, conditioned to hit $A$, first enters $A$; $\harm_A^\infty$ is its limit as $\dist(x,A;G) \to \infty$, when this limit exists. For the plane data $\cK = \{J_x:x \in \cV(G)\}$ of Section~\ref{subsec:annular-intersections} and finite connected $A \subset \cV(G)$, set $U_A^{\cK} := \bigcup_{x \in A}J_x$. Note that when the plane data $\cK$ is a cell configuration $\cH$, $U_A^{\cK}$ is connected. Set $D(A) := \diam(U_A^{\cH})$.

\paragraph{Acknowledgements} The author would like to thank Ewain Gwynne for helpful conversations throughout the duration of this project, and Leonardo Bonanno for preparing the associated formalization of this paper.

\paragraph{AI Disclosure} ChatGPT 5.6 and 6 Pro were used to expand a detailed proof sketch into an initial draft early in this project, for subsequent proof revisions, to generate Figures~\ref{fig:good-scales} and \ref{fig:polygonal-path}, and for the editing of prose, proofs and \LaTeX. All results have been verified by the author, and any errors herein are the sole responsibility of the author. A formalization of this paper prepared by Leonardo Bonanno is available in~\cite{bonannoFormalizations}.

\section{Preliminaries}\label{sec:prelims}
Now, we discuss in greater detail the relationship between the mated CRT map and SLE/LQG, after which we present some preliminary results on the mated CRT map and ergodic scale-free random environments. Before that, we define SLE, the Gaussian free field (GFF) and LQG.

\subsection{Basics of LQG and SLE}
Now, we describe the basics of certain variants of SLE, GFF and LQG. We cover only the material necessary to describe the canonical embedding and cellular embedding with which the mated CRT map comes equipped. Since we rely on many established properties of the mated CRT map, but do not exploit the relationship between LQG/SLE and the mated CRT map very directly, we keep our presentation brief, and only state those definitions and properties necessary to interpret the relevant mated CRT properties. Further references on SLE and LQG are given. In what follows, fix $\gamma \in (0,2)$ and set $Q := 2/\gamma + \gamma/2$ and $\kappa := 16/\gamma^2$.

\subsubsection{Space filling SLE}
SLE is a family of random fractal curves that arise as scaling limits of interfaces in critical lattice models \cite{schramm0, smirnov-cardy, lawler2011conformal, smirnov2010conformal} and are indexed by a parameter $\kappa > 0$, which varies the ``roughness'' of the curve. We give a very brief introduction to these random curves here, and refer to \cite{lawler2008conformally} for further details.

In this paper, we deal exclusively with a variant of SLE known as whole-plane space-filling $\mathrm{SLE}_\kappa$ from $\infty$ to $\infty$ with $\kappa = 16/\gamma^2$ when $\gamma \in (0, 2)$ (see e.g. \cite[Section~1.4.1]{duplantier2021liouville} for the original definition in the whole-plane case, or \cite[Section~2.1.3]{gwynne2019harmonicfunctionsmatedcrtmaps} for a more recent discussion). This is a random continuous curve $\eta: \R \to \C$ with $\eta(\R) = \C$ and $|\eta(t)| \to \infty$ as $|t| \to \infty$. The space-filling construction is identical to ordinary whole plane SLE when $\kappa \geq 8$, and is a slight modification when $4 < \kappa < 8$. Its law is invariant under translations, rotations, and positive dilations.

\subsubsection{Whole plane GFF and LQG}
Now, we introduce the GFF and LQG. Let $\mathcal{H}(\C)$ be the completion of $\{f \in C_c^\infty(\C): \int_\C f(z)\,dz = 0\}$ for the Dirichlet inner product $(f,g)_\nabla := (2 \pi)^{-1} \int_\C \nabla f \cdot \nabla g\,dz$. Given an orthonormal basis $(f_j)$ and independent standard Gaussian variables $(\alpha_j)$, the formal sum
\[
    h := \sum_{j \geq 1} \alpha_jf_j
\]
defines the whole-plane GFF, a random distribution modulo additive constants. For each $f \in \mathcal{H}(\C)$, the pairing $(h,f)_\nabla := \sum_{j \geq 1} \alpha_j(f_j,f)_\nabla$ converges in $L^2$, and $\E[(h,f)_\nabla(h,g)_\nabla] = (f,g)_\nabla$; see \cite[Section~2.1.1]{gwynne2019harmonicfunctionsmatedcrtmaps}.

For $z \in \C$ and $r > 0$, let $h_r(z)$ denote the average of $h$ with respect to the uniform probability measure on $\partial B_r(z)$, defined in \cite[Section~3.1]{duplantier2011liouville}. We fix the additive constant by normalizing $h_1(0) = 0$.

For a conformal map $f:D \to \widetilde D$, define $f \bullet_\gamma h := h \circ f^{-1} + Q \log|(f^{-1})'|$. Write $(D,h) \sim_\gamma (\widetilde D, \widetilde h)$ if $\widetilde h = f \bullet_\gamma h$ for such an $f$. A \emph{$\gamma$-LQG surface} is an equivalence class for this relation, and a representative is called an \emph{embedding}. For marked surfaces, $f$ must send each marked point to its corresponding marked point; see e.g. \cite[Section~2.3]{ang2026quantum}.

The LQG area measure is formally
$$
    \cA_h^\gamma := \lim_{\epsilon \to 0} \epsilon^{\gamma^2/2}e^{\gamma h_\epsilon(z)}\,dz,
$$
where the convergence is a.s.\ in the local weak topology on measures on $D$ along dyadic $\epsilon$. It satisfies $\cA_{f \bullet_\gamma h}^\gamma = f_*\cA_h^\gamma$ and $\cA_{h + c}^\gamma = e^{\gamma c} \cA_h^\gamma$, and is therefore well defined for a $\gamma$-LQG surface (see \cite[Sections 2 and 3]{duplantier2011liouville}). We use $\operatorname{Area}$ for the Lebesgue measure.

To study the mated CRT map, we will need a variant of LQG surface called an $\alpha$-quantum cone with $\alpha \in \{0, \gamma \}$. Let $B$ be standard one-dimensional Brownian motion and, independently of $B$, let $\widehat{B}$ have the law of Brownian motion conditioned so that $\widehat{B}_t + (Q - \alpha)t > 0$ for all $t > 0$. Let
$$
    A_t :=
    \begin{cases}
        B_t + \alpha t, & t \geq 0,\\
        \widehat{B}_{-t} + \alpha t, & t < 0.
    \end{cases}
$$
The circle-average embedding of the $\alpha$-quantum cone is defined as the field $h$ for which $t \mapsto h_{e^{-t}}(0)$ has the law of $A_t$, while $h - h_{|\cdot|}(0)$ is independent of $A_t$ and agrees in law with the corresponding projection of a whole-plane GFF. This field satisfies $\sup \{r > 0: h_r(0) + Q \log r = 0\} = 1$ (see \cite[Definition~4.10]{duplantier2021liouville} and \cite[Section~2.3.1]{gwynne2021tutteembeddingmatedcrtmap}). In this paper, we will take $\alpha \in \{0, \gamma \}$.

\subsection{The mated CRT and LQG/SLE}\label{sec:crt-lqg}
Let $h$ be the circle-average embedding of the $\gamma$-quantum cone defined above, and $\eta$ an independent whole-plane space-filling $\mathrm{SLE}_\kappa$. We parametrize $\eta$ by $\gamma$-LQG area: $\eta(0) = 0$ and $\cA_h^{\gamma}(\eta([s,t])) = t - s$ for $s < t$. This parametrization is the so-called \textit{quantum natural time} parametrization of SLE, and has a number of useful properties (see e.g. \cite[Section 2.5]{ang2026quantum}). For $\epsilon > 0$ define
\begin{equation*}
  H_x^\epsilon := \eta([x - \epsilon,x]),
  \qquad x \in \epsilon \Z.
\end{equation*}
These are compact connected sets of measure $\epsilon$ under $\cA_h^{\gamma}$. The structure graph of the cells $\{H_x^\epsilon \}$ is defined as follows: two distinct vertices are adjacent when their cells share a nontrivial boundary arc, consecutive indices are joined by one edge, nonconsecutive indices are joined by two edges when both left and right boundary adjacency occurs \cite[Section 1.2]{gwynne2019harmonicfunctionsmatedcrtmaps}. The resulting graph has the law of the whole-plane $\gamma$-mated-CRT map $\mathcal{G}^\epsilon$ (see e.g. \cite[Proposition 2.2]{gwynne2021tutteembeddingmatedcrtmap}, \cite[Section 1.2]{gwynne2019harmonicfunctionsmatedcrtmaps}). This observation allows one to study the $\gamma$-mated CRT map via the $\gamma$-quantum cone and space filling SLE, and to couple the $\gamma$-mated CRT maps for different $\epsilon$.

For $D \subset \C$, write $\cV_{\mathcal{G}^\epsilon}(D) := \{x \in \epsilon \Z:H_x^\epsilon \cap D \neq \varnothing \}$. At unit mesh, the root is $o = 0$; set $H_x := H_x^1$, $\mathcal{H} := \{H_x:x \in \Z \}$, and $U_A^{\mathcal{H}} := \bigcup_{x \in A}H_x$. Distinct cells overlap only on sets of zero quantum area, so $\cA_h^\gamma(U_A^{\mathcal{H}}) = |A|$ for finite $A \subset \Z$. This is the cellular embedding used in the verification of Assumption~\ref{assumption:mesoscopic}.

\subsection{Properties of the mated CRT map}
We retain the SLE/LQG realization and notation of the preceding subsection and record the properties needed for the mesoscopic argument. We start with control of the random walk and harmonic measure.

\begin{proposition}[{\normalfont(\cite[Theorem~1.4]{gwynne2019harmonicfunctionsmatedcrtmaps})}]
\label{prop:mated-recurrence}
For every $\epsilon > 0$, the simple random walk on $\mathcal{G}^\epsilon$ is recurrent almost surely.
\end{proposition}

We denote the total variation (TV) distance by $\| \cdot \|_{\mathrm{TV}}$.

\begin{proposition}[{\normalfont(\cite[Theorem~1.1 and Section~9]{berestycki2024harnack}, \cite[Proposition~1.7]{gwynne2021external})}]\label{prop:mated-hm}
Almost surely, for every finite nonempty connected $A \subset \Z$ we have
\begin{equation*}
  \lim_{m \to \infty} \sup_{\substack{w \in \cV(\mathcal{G}^1) \\ \dist(w,A; \mathcal{G}^1) \geq m}} \| \harm^w_A - \harm^\infty_A \|_{\mathrm{TV}} = 0
\end{equation*}
where $\harm^\infty_A$ is a unique measure on $\bigcup_{v \in A} \cE_A(v)$.
\end{proposition}

For $b > 0$, define $R_b := \sup \{r > 0:h_r(0) + Q \log r = \gamma^{-1} \log b\}$, $h^b(\cdot) := h(R_b \cdot) + Q \log R_b - \gamma^{-1} \log b$, and $\eta^b(t) := R_b^{-1} \eta(bt)$.

\begin{proposition}[{\normalfont(\cite[Proposition~4.13(i)]{duplantier2021liouville})}]\label{prop:quantum-cone-scaling}
For every $b > 0$ we have
\begin{equation*}
  (h^b, \eta^b) \overset d = (h, \eta).
\end{equation*}
As a result, the $1/b$-cell configuration for $(h^b, \eta^b)$ is, after multiplication by $R_b$, the $1$-cell configuration for $(h, \eta)$.
\end{proposition}
Define
\begin{equation*}
  a_\gamma := \frac{1}{\gamma(Q - \gamma)} = \frac{1}{2 - \gamma^2/2}.
\end{equation*}
The following result from \cite{gwynne2020anomalous} controls the asymptotics of $R_b$.

\begin{proposition}[{\normalfont(\cite[Lemma~2.6 and the proof of Theorem~1.4]{gwynne2020anomalous})}]\label{prop:radius-asymp}
Almost surely,
\begin{equation*}
  \lim_{k \to \infty} \frac{\log R_{2^k}}{k \log 2} = a_\gamma.
\end{equation*}
\end{proposition}
Next, we need control over the number of vertices in $B_{R_b}(0)$ with $b = 2^k$.
\begin{proposition}[{\normalfont(\cite[proof of Theorem~1.4]{gwynne2020anomalous};
\cite[Lemma~A.3]{gwynne2019distance})}]
\label{prop:mated-euclidean-cell-count}
For every $\delta > 0$, almost surely
\begin{equation*}
  |\cV_{\mathcal{G}^1}(B_{R_{2^k}}(0))|
  \leq 2^{(1 + \delta)k}
\end{equation*}
for all sufficiently large $k$.
\end{proposition}

The above estimate appears immediately after Equation (3.31) in \cite[]{gwynne2020anomalous} with probability at least $1 - O(2^{-\alpha k})$ for some $\alpha > 0$, and Borel-Cantelli gives the a.s.\ claim.

\subsection{Mesoscopic control of random walk for mated CRT}
The following estimate is extracted from \cite[Proposition 3.6, Lemma 2.7]{gwynne2019harmonicfunctionsmatedcrtmaps} and the proof of \cite[Lemma 3.11]{gwynne2019harmonicfunctionsmatedcrtmaps}.

\begin{proposition}\label{prop:mated-crt-mesoscopic}
There exist constants
\begin{equation*}
  \rho \in (0,1), \qquad p, \alpha, \beta, \chi,q > 0, \qquad \beta < \min \{ \chi,q\}, \qquad q < \frac{2}{(2 + \gamma)^2}
\end{equation*}
such that the following holds. Define
\begin{equation}
  \mathcal{Z}_\epsilon := (\epsilon^\chi \Z^2) \cap B_\rho(0),
  \qquad
  \mathcal{R}_\epsilon := (\epsilon^\chi \Z) \cap [\epsilon^\beta, \rho].
 \label{eq:input-grids}
\end{equation}
The following hold simultaneously with probability at least $1 - O_\epsilon(\epsilon^\alpha)$:
\begin{enumerate}
\item Every cell intersecting $B_\rho(0)$ has Euclidean diameter at most $\epsilon^q$.
\item For every $(z,r) \in \mathcal{Z}_\epsilon \times \mathcal{R}_\epsilon$ with $B_{2r}(z) \subset B_\rho(0)$, and for every $x \in \cV_{\mathcal{G}^\epsilon}(\partial B_r(z))$, let $X_n$ be simple random walk started from $x$ and $\tau$ its first exit time from $\cV_{\mathcal{G}^\epsilon}(\A_{r/2, 3r/2}(z))$. With quenched probability at least $p$ we have:
\begin{enumerate}
\item The random walk disconnects $\cV_{\mathcal{G}^\epsilon}(\partial B_{3r/4}(z))$ from $\cV_{\mathcal{G}^\epsilon}(\partial B_{5r/4}(z))$ before $\tau$.
\item Every vertex visited before $\tau$ belongs to $\cV_{\mathcal{G}^\epsilon}(\A_{r/2, 3r/2}(z))$.
\item The union of cells visited up to and including $\tau$ is contained in $\A_{r/3, 5r/3}(z)$.
\end{enumerate}
\end{enumerate}
\end{proposition}

The proof of \cite[Lemma 3.11]{gwynne2019harmonicfunctionsmatedcrtmaps} applies \cite[Proposition 3.6]{gwynne2019harmonicfunctionsmatedcrtmaps} to obtain $(2a)$ and $(2b)$ with quenched success probability at least $p > 0$. The final union bound makes these estimates simultaneous over the grids in \eqref{eq:input-grids}: choosing $\chi > 0$ sufficiently small and $\beta$ below both $\chi$ and the lower-scale exponent in that lemma gives a failure probability at most $O_\epsilon(\epsilon^\alpha)$. \cite[Lemma 2.7]{gwynne2019harmonicfunctionsmatedcrtmaps} gives $(1)$, after choosing $\beta < q < 2/(2 + \gamma)^2$ and decreasing $\alpha$ if necessary. Since every admissible radius satisfies $r \geq \epsilon^\beta$, one has $\epsilon^q = o(r)$. So, $(1)$ and $(2b)$ imply that the union of cells visited up to and including $\tau$ is contained in $\A_{r/3, 5r/3}(z)$ for all sufficiently small $\epsilon$, which establishes $(2c)$.

\subsection{Ergodic scale-free random environments}
Lastly, we introduce the invariance principle for ergodic scale-free random environments. First, we need a definition.

\begin{definition}\label{def:conductance-weighted-walk}
For a good embedded lattice $\mathcal{M}$, the conductance weighted random walk on $\mathcal{M}^*$ has transition probabilities $P(H,H') = \mathfrak c(H,H')/\pi(H)$. Here $\mathcal{M}^*$ is viewed as a weighted graph with vertices given by the faces of $\mathcal{M}$, with one edge of conductance $\mathfrak c(H,H')$ for each adjacent pair of faces. This Markov chain on the faces of $\mathcal{M}$ is the random walk defined in \cite[Theorem 3.10]{gwynne2021invarianceprincipleergodicscalefree}.
\end{definition}

The next result we recall is a main result from \cite{gwynne2021invarianceprincipleergodicscalefree}; it is the invariance principle for ergodic-scale free random environments. In particular, this proposition combines Theorems 1.5 and 3.10 from \cite{gwynne2021invarianceprincipleergodicscalefree}.

\begin{proposition}\label{prop:ergodic}
Let $\mathcal{M}$ be a good embedded lattice (see Definition~\ref{def:ergodic}), and choose a point $p_H \in \operatorname{int}(H)$ from each face, with $p_{H_0} = 0$. Let $(X_n)$ be the conductance weighted random walk of $\mathcal{M}^*$ in Definition~\ref{def:conductance-weighted-walk}. For a starting face $H$, let $\widehat{Y}^{H;R,L}$ be the piecewise-linear interpolation of $n \mapsto p_{X_n}$, started from $p_H$ and stopped upon first leaving $B_{LR}(0)$. There exists a positive-definite covariance matrix $\Sigma$ such that, a.s.\ in $\mathcal{M}$, for every $0 < C < L < \infty$, we have
$$
  \sup_{H:\,|p_H| \leq CR}d_{\textrm{Pr}}(\cL_H(R^{-1} \widehat{Y}^{H;R,L} | \, \mathcal{M}), \cL(B_\Sigma^{p_H/R;L}))
  \to 0 \text{ as } R \to \infty
$$
where $B_\Sigma^{z;L}$ is Brownian motion with covariance matrix $\Sigma$ started from $z$ and stopped upon first leaving $B_L(0)$, and $d_{\textrm{Pr}}$ is induced by $1 \wedge d_F$. Convergence is again with respect to the topology on curves modulo time parametrization.
\end{proposition}

Now we proceed to the proofs of the blackbox results.

\section{Control of diameter via annular intersection probabilities}\label{sec:arg-template}
In this section, we demonstrate how to convert annular intersection probabilities into harmonic measure bounds, and harmonic measure bounds into diameter bounds. Sections~\ref{sec:pf-sublinear} and \ref{sec:pf-poly-sublinear} will invoke these intermediary results to establish the main theorems. The section concludes with a brief subsection describing the change of initial vertex of the DLA growth.

\subsection{Annular intersection probabilities to harmonic measure bounds}\label{subsec:annular-intersections}
The results of this section are stated at a fixed scale. In the next subsection, these results will be applied at an arbitrary scale, yielding the large $n$ conclusions.

Assumptions~\ref{assumption:rw-conv} and ~\ref{assumption:mesoscopic} pertain to slightly different embedding data for the planar map (a proper planar embedding versus a cellular embedding), and the arguments obey this structure. To reconcile these differences, we define the \textit{plane data} as follows: let $\cK = \{J_x: x \in \cV(G)\}$ be a locally finite collection of nonempty bounded connected subsets of $\C$, indexed by $\cV(G)$, such that $J_x \cap J_y \neq \emptyset$ if $x \sim y$. In Section~\ref{sec:pf-sublinear}, $J_x$ will be the union of $\{ \Phi(x)\}$ and the interiors of the embedded edges incident to $x$, and in Section~\ref{sec:pf-poly-sublinear} we choose $\cK$ as the cellular embedding. For $D \subset \C$, define
\begin{equation*}
    \cV_{\cK}(D) := \{x \in \cV(G): J_x \cap D \neq \emptyset \}
\end{equation*}
as the vertices with planar data intersecting a given Euclidean subset.

Now, we state precisely the multiscale structure, at a fixed non-normalized scale, which will yield asymptotically shrinking harmonic measure bounds.

\begin{definition}[Good scales]\label{def:good-scales}
    Fix $K \in \N$ and $p \in (0, 1)$, and let $A \subset \cV(G)$ be finite and connected. Choose $v \in A$. $(K,p)$-good scales for $(A,v)$ consist of radii 
    \begin{equation}\label{eq:good-radii}
        r_1 > r_2 > \cdots > r_K > 0, \qquad r_{j + 1} \leq \frac{r_j}{8},
    \end{equation}
    and a center $z \in \C$ such that:
    \begin{enumerate}
        \item For all $J_x$ such that $J_x \cap B_{2r_1}(z) \neq \emptyset$ we have $\diam(J_x) < r_K/100$.
        \item $\cup_{x \in N(v)} \, J_x \subset B_{r_K/20}(z).$
        \item For all $j \in [K]$, there exists a graph path in $A$ joining $\cV_{\cK}(\partial B_{3r_j/4}(z))$ to $\cV_{\cK}(\partial B_{5r_j/4}(z))$.
        \item For all $j \in [K]$ and all $x \in \cV_{\cK}(\partial B_{r_j}(z))$, a simple random walk $X_n$ started from $x$ and run until its first exit time $\tau$ from $W_j := \cV_{\cK}(\A_{r_j/2, 3r_j/2}(z))$ has probability at least $p$ that it disconnects $\cV_{\cK}(\partial B_{3r_j/4}(z))$ from $\cV_{\cK}(\partial B_{5r_j/4}(z))$ before $\tau$ and $\cup_{m = 0}^{\tau} \, J_{X_m} \subset \A_{r_j/3, 5r_j/3}(z)$.
    \end{enumerate}
\end{definition}
Note that the above definition admits at most \textit{logarithmically many} good scales relative to the ratio $r_1/r_K$. In principle one could achieve $K \gg \log(r_1/r_K)$ without \eqref{eq:good-radii}, but we do not encounter this scenario in this paper and so work with the slightly more restrictive definition.

\begin{figure}[ht]
\centering
\begin{tikzpicture}[
    scale=1.25,
    line cap=round,
    line join=round,
    map edge/.style={draw=black!40,line width=.4pt},
    map vertex/.style={circle,fill=black!65,inner sep=.65pt},
    cluster/.style={draw=black,line width=1.1pt},
    vertex/.style={circle,fill=black,inner sep=1.05pt},
    label/.style={font=\scriptsize,inner sep=1.2pt},
    outer annulus/.style={red,densely dotted,line width=.75pt},
    inner annulus/.style={blue,densely dotted,line width=.75pt}]

    \coordinate (n20) at (-1.466,-2.957);
    \coordinate (n30) at (-0.633,-2.785);
    \coordinate (n40) at (-0.041,-2.933);
    \coordinate (n50) at (0.719,-2.881);
    \coordinate (n60) at (1.533,-2.889);
    \coordinate (n11) at (-2.211,-2.145);
    \coordinate (n21) at (-1.39,-2.075);
    \coordinate (n31) at (-0.821,-2.118);
    \coordinate (n41) at (-0.105,-2.263);
    \coordinate (n51) at (0.811,-2.097);
    \coordinate (n61) at (1.357,-2.156);
    \coordinate (n71) at (2.19,-2.199);
    \coordinate (n02) at (-2.87,-1.417);
    \coordinate (n12) at (-2.12,-1.492);
    \coordinate (n22) at (-1.497,-1.388);
    \coordinate (n32) at (-0.796,-1.397);
    \coordinate (n42) at (-0.1,-1.372);
    \coordinate (n52) at (0.813,-1.457);
    \coordinate (n62) at (1.359,-1.392);
    \coordinate (n72) at (2.187,-1.407);
    \coordinate (n82) at (2.944,-1.43);
    \coordinate (n03) at (-3.12,-0.74);
    \coordinate (n13) at (-2.154,-0.807);
    \coordinate (n23) at (-1.404,-0.637);
    \coordinate (n33) at (-0.79,-0.764);
    \coordinate (n43) at (0.1,-0.656);
    \coordinate (n53) at (0.765,-0.784);
    \coordinate (n63) at (1.436,-0.621);
    \coordinate (n73) at (2.145,-0.654);
    \coordinate (n83) at (2.873,-0.829);
    \coordinate (n04) at (-2.934,-0.039);
    \coordinate (n14) at (-2.093,0.056);
    \coordinate (n24) at (-1.403,-0.103);
    \coordinate (n34) at (-0.657,0.062);
    \coordinate (n44) at (0,0);
    \coordinate (n54) at (0.776,0.044);
    \coordinate (n64) at (1.455,-0.097);
    \coordinate (n74) at (2.092,-0.086);
    \coordinate (n84) at (2.874,-0.102);
    \coordinate (n05) at (-2.985,0.691);
    \coordinate (n15) at (-2.235,0.788);
    \coordinate (n25) at (-1.33,0.686);
    \coordinate (n35) at (-0.666,0.799);
    \coordinate (n45) at (0.04,0.776);
    \coordinate (n55) at (0.729,0.664);
    \coordinate (n65) at (1.502,0.67);
    \coordinate (n75) at (2.179,0.719);
    \coordinate (n85) at (2.822,0.727);
    \coordinate (n06) at (-2.79,1.465);
    \coordinate (n16) at (-2.218,1.371);
    \coordinate (n26) at (-1.35,1.514);
    \coordinate (n36) at (-0.737,1.474);
    \coordinate (n46) at (0.109,1.537);
    \coordinate (n56) at (0.691,1.394);
    \coordinate (n66) at (1.454,1.506);
    \coordinate (n76) at (2.228,1.452);
    \coordinate (n86) at (2.824,1.43);
    \coordinate (n17) at (-2.081,2.191);
    \coordinate (n27) at (-1.358,2.067);
    \coordinate (n37) at (-0.646,2.135);
    \coordinate (n47) at (-0.015,2.21);
    \coordinate (n57) at (0.685,2.269);
    \coordinate (n67) at (1.516,2.175);
    \coordinate (n77) at (2.164,2.257);
    \coordinate (n28) at (-1.432,2.813);
    \coordinate (n38) at (-0.63,2.893);
    \coordinate (n48) at (-0.018,2.907);
    \coordinate (n58) at (0.692,2.956);
    \coordinate (n68) at (1.437,2.957);

    \foreach \rin/\rout in {0.60/0.90,1.42/1.73,2.34/2.65} {
        \path[fill=black!6,even odd rule]
            (n44) circle[radius=\rout] (n44) circle[radius=\rin];
    }

    \draw[map edge] (n20) -- (n30);
    \draw[map edge] (n20) -- (n21);
    \draw[map edge] (n30) -- (n40);
    \draw[map edge] (n30) -- (n31);
    \draw[map edge] (n40) -- (n50);
    \draw[map edge] (n40) -- (n41);
    \draw[map edge] (n50) -- (n60);
    \draw[map edge] (n50) -- (n51);
    \draw[map edge] (n60) -- (n61);
    \draw[map edge] (n11) -- (n21);
    \draw[map edge] (n11) -- (n12);
    \draw[map edge] (n21) -- (n31);
    \draw[map edge] (n21) -- (n22);
    \draw[map edge] (n31) -- (n41);
    \draw[map edge] (n31) -- (n32);
    \draw[map edge] (n41) -- (n51);
    \draw[map edge] (n41) -- (n42);
    \draw[map edge] (n51) -- (n61);
    \draw[map edge] (n51) -- (n52);
    \draw[map edge] (n61) -- (n71);
    \draw[map edge] (n61) -- (n62);
    \draw[map edge] (n71) -- (n72);
    \draw[map edge] (n02) -- (n12);
    \draw[map edge] (n02) -- (n03);
    \draw[map edge] (n12) -- (n22);
    \draw[map edge] (n12) -- (n13);
    \draw[map edge] (n22) -- (n32);
    \draw[map edge] (n22) -- (n23);
    \draw[map edge] (n32) -- (n42);
    \draw[map edge] (n32) -- (n33);
    \draw[map edge] (n42) -- (n52);
    \draw[map edge] (n42) -- (n43);
    \draw[map edge] (n52) -- (n62);
    \draw[map edge] (n52) -- (n53);
    \draw[map edge] (n62) -- (n72);
    \draw[map edge] (n62) -- (n63);
    \draw[map edge] (n72) -- (n82);
    \draw[map edge] (n72) -- (n73);
    \draw[map edge] (n82) -- (n83);
    \draw[map edge] (n03) -- (n13);
    \draw[map edge] (n03) -- (n04);
    \draw[map edge] (n13) -- (n23);
    \draw[map edge] (n13) -- (n14);
    \draw[map edge] (n23) -- (n33);
    \draw[map edge] (n23) -- (n24);
    \draw[map edge] (n33) -- (n43);
    \draw[map edge] (n33) -- (n34);
    \draw[map edge] (n43) -- (n53);
    \draw[map edge] (n43) -- (n44);
    \draw[map edge] (n53) -- (n63);
    \draw[map edge] (n53) -- (n54);
    \draw[map edge] (n63) -- (n73);
    \draw[map edge] (n63) -- (n64);
    \draw[map edge] (n73) -- (n83);
    \draw[map edge] (n73) -- (n74);
    \draw[map edge] (n83) -- (n84);
    \draw[map edge] (n04) -- (n14);
    \draw[map edge] (n04) -- (n05);
    \draw[map edge] (n14) -- (n24);
    \draw[map edge] (n14) -- (n15);
    \draw[map edge] (n24) -- (n34);
    \draw[map edge] (n24) -- (n25);
    \draw[map edge] (n34) -- (n44);
    \draw[map edge] (n34) -- (n35);
    \draw[map edge] (n44) -- (n54);
    \draw[map edge] (n44) -- (n45);
    \draw[map edge] (n54) -- (n64);
    \draw[map edge] (n54) -- (n55);
    \draw[map edge] (n64) -- (n74);
    \draw[map edge] (n64) -- (n65);
    \draw[map edge] (n74) -- (n84);
    \draw[map edge] (n74) -- (n75);
    \draw[map edge] (n84) -- (n85);
    \draw[map edge] (n05) -- (n15);
    \draw[map edge] (n05) -- (n06);
    \draw[map edge] (n15) -- (n25);
    \draw[map edge] (n15) -- (n16);
    \draw[map edge] (n25) -- (n35);
    \draw[map edge] (n25) -- (n26);
    \draw[map edge] (n35) -- (n45);
    \draw[map edge] (n35) -- (n36);
    \draw[map edge] (n45) -- (n55);
    \draw[map edge] (n45) -- (n46);
    \draw[map edge] (n55) -- (n65);
    \draw[map edge] (n55) -- (n56);
    \draw[map edge] (n65) -- (n75);
    \draw[map edge] (n65) -- (n66);
    \draw[map edge] (n75) -- (n85);
    \draw[map edge] (n75) -- (n76);
    \draw[map edge] (n85) -- (n86);
    \draw[map edge] (n06) -- (n16);
    \draw[map edge] (n16) -- (n26);
    \draw[map edge] (n16) -- (n17);
    \draw[map edge] (n26) -- (n36);
    \draw[map edge] (n26) -- (n27);
    \draw[map edge] (n36) -- (n46);
    \draw[map edge] (n36) -- (n37);
    \draw[map edge] (n46) -- (n56);
    \draw[map edge] (n46) -- (n47);
    \draw[map edge] (n56) -- (n66);
    \draw[map edge] (n56) -- (n57);
    \draw[map edge] (n66) -- (n76);
    \draw[map edge] (n66) -- (n67);
    \draw[map edge] (n76) -- (n86);
    \draw[map edge] (n76) -- (n77);
    \draw[map edge] (n17) -- (n27);
    \draw[map edge] (n27) -- (n37);
    \draw[map edge] (n27) -- (n28);
    \draw[map edge] (n37) -- (n47);
    \draw[map edge] (n37) -- (n38);
    \draw[map edge] (n47) -- (n57);
    \draw[map edge] (n47) -- (n48);
    \draw[map edge] (n57) -- (n67);
    \draw[map edge] (n57) -- (n58);
    \draw[map edge] (n67) -- (n77);
    \draw[map edge] (n67) -- (n68);
    \draw[map edge] (n28) -- (n38);
    \draw[map edge] (n38) -- (n48);
    \draw[map edge] (n48) -- (n58);
    \draw[map edge] (n58) -- (n68);
    \draw[map edge] (n20) -- (n31);
    \draw[map edge] (n40) -- (n31);
    \draw[map edge] (n50) -- (n41);
    \draw[map edge] (n50) -- (n61);
    \draw[map edge] (n21) -- (n12);
    \draw[map edge] (n31) -- (n22);
    \draw[map edge] (n41) -- (n32);
    \draw[map edge] (n61) -- (n52);
    \draw[map edge] (n61) -- (n72);
    \draw[map edge] (n12) -- (n03);
    \draw[map edge] (n12) -- (n23);
    \draw[map edge] (n32) -- (n23);
    \draw[map edge] (n32) -- (n43);
    \draw[map edge] (n62) -- (n53);
    \draw[map edge] (n72) -- (n83);
    \draw[map edge] (n13) -- (n04);
    \draw[map edge] (n13) -- (n24);
    \draw[map edge] (n23) -- (n34);
    \draw[map edge] (n33) -- (n44);
    \draw[map edge] (n73) -- (n64);
    \draw[map edge] (n04) -- (n15);
    \draw[map edge] (n24) -- (n15);
    \draw[map edge] (n34) -- (n25);
    \draw[map edge] (n34) -- (n45);
    \draw[map edge] (n54) -- (n65);
    \draw[map edge] (n64) -- (n75);
    \draw[map edge] (n74) -- (n85);
    \draw[map edge] (n05) -- (n16);
    \draw[map edge] (n35) -- (n26);
    \draw[map edge] (n35) -- (n46);
    \draw[map edge] (n45) -- (n56);
    \draw[map edge] (n55) -- (n66);
    \draw[map edge] (n75) -- (n66);
    \draw[map edge] (n75) -- (n86);
    \draw[map edge] (n26) -- (n17);
    \draw[map edge] (n26) -- (n37);
    \draw[map edge] (n56) -- (n47);
    \draw[map edge] (n56) -- (n67);
    \draw[map edge] (n66) -- (n77);
    \draw[map edge] (n37) -- (n28);
    \draw[map edge] (n37) -- (n48);
    \draw[map edge] (n57) -- (n48);
    \draw[map edge] (n67) -- (n58);

    \foreach \rin/\rout in {0.60/0.90,1.42/1.73,2.34/2.65} {
        \draw[inner annulus] (n44) circle[radius=\rin];
        \draw[outer annulus] (n44) circle[radius=\rout];
    }

    \draw[cluster] (n34) -- (n24) -- (n13) -- (n03);
    \draw[cluster] (n13) -- (n12) -- (n11);
    \draw[cluster] (n34) -- (n25) -- (n26) -- (n17);
    \draw[cluster] (n26) -- (n37) -- (n38);
    \draw[cluster] (n34) -- (n44) -- (n54) -- (n64) -- (n74);
    \draw[cluster] (n54) -- (n55) -- (n56) -- (n57);
    \draw[cluster] (n34) -- (n33) -- (n32) -- (n31) -- (n30);
    \draw[cluster] (n32) -- (n42) -- (n52) -- (n61);

    \foreach \p in {n20,n30,n40,n50,n60,n11,n21,n31,n41,n51,n61,n71,n02,n12,n22,n32,n42,n52,n62,n72,n82,n03,n13,n23,n33,n43,n53,n63,n73,n83,n04,n14,n24,n34,n44,n54,n64,n74,n84,n05,n15,n25,n35,n45,n55,n65,n75,n85,n06,n16,n26,n36,n46,n56,n66,n76,n86,n17,n27,n37,n47,n57,n67,n77,n28,n38,n48,n58,n68}
        \node[map vertex] at (\p) {};
    \foreach \p in {n03,n11,n12,n13,n17,n24,n25,n26,n30,n31,n32,n33,n34,n37,n38,n42,n44,n52,n54,n55,n56,n57,n61,n64,n74}
        \node[vertex] at (\p) {};

    \node[label,anchor=south,yshift=4pt] at (n34) {$o$};
    \node[label,anchor=north west,xshift=2pt,yshift=-2pt] at (n44) {$v$};
    \node[label] at (-1.62,1.10) {$A$};

    \node[label] at (1.2475,2.1607) {$\mathbb{A}_1$};
    \node[label] at (1.22,.995) {$\mathbb{A}_2$};
    \node[label] at (.575,-.482) {$\mathbb{A}_3$};
\end{tikzpicture}

\caption{$(K, p)$-good scales for the pair $(A, v)$ with $K = 3$. The good scales consist of $3$ disjoint annuli $\mathbb{A}_1, \mathbb{A}_2$ and $\mathbb{A}_3$, with inner and outer boundaries colored blue and red respectively. A segment of the ambient graph $G$ is drawn with gray vertices and edges, and a segment of $A$ is drawn with black vertices and edges.}
\label{fig:good-scales}
\end{figure}
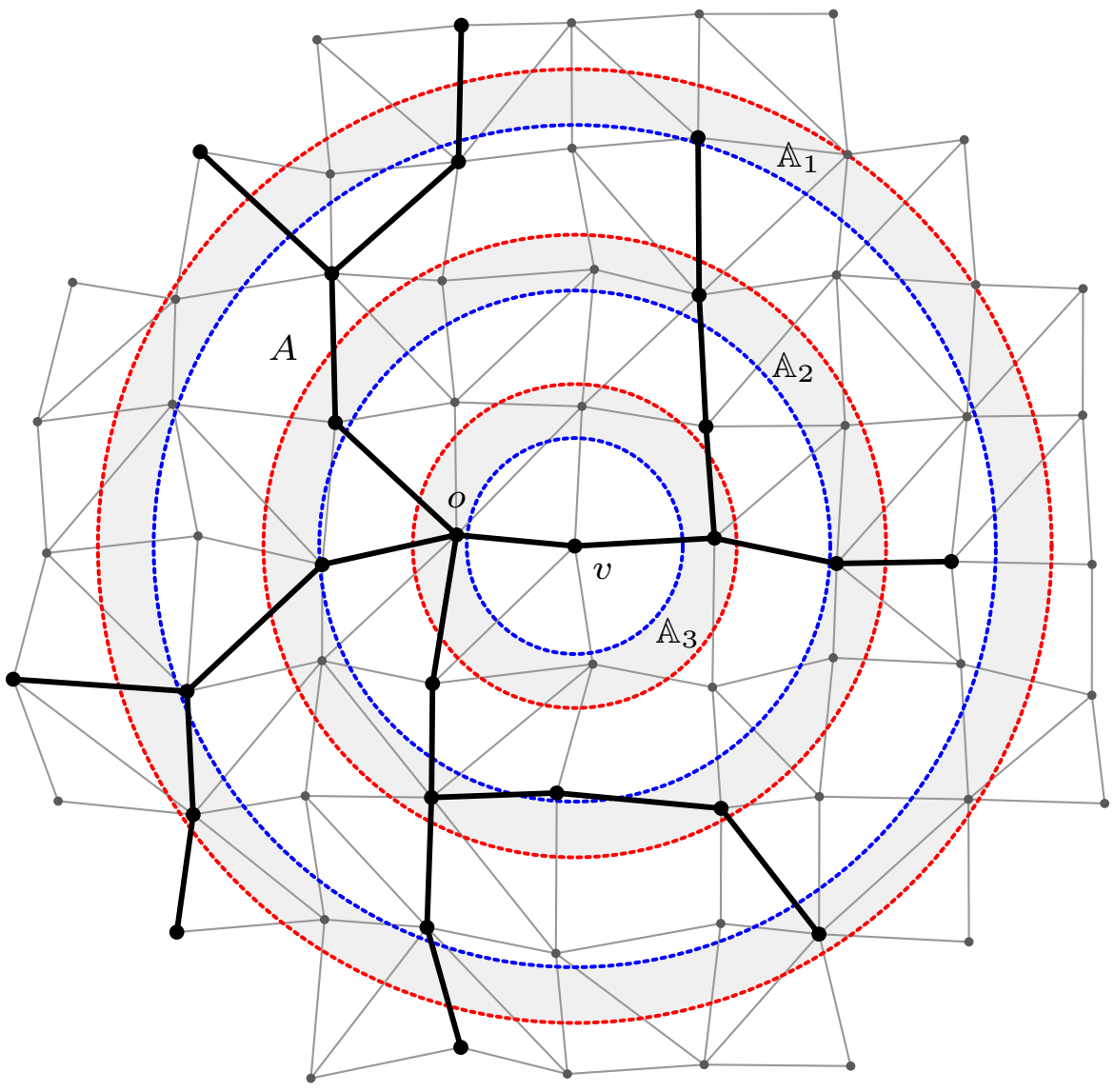

As we organize scales in a Euclidean fashion, i.e. using Euclidean annuli, the bound on diameters of $J_{x}$ establishes a compatibility condition of the plane data with the good scales, and ensures that a random walk from $\infty$ must pass through each individual scale successively. Each scale then provides with probability at least $p$ the possibility of an annular disconnection, and combining with the Markov property of the random walk yields the final harmonic measure bound. This is illustrated in Figure~\ref{fig:good-scales}, and is the content of the next proposition.

\begin{proposition}\label{prop:helper-harm}
Suppose there exist $(K,p)$-good scales for $(A,v)$. Then for all $x \in \cV(G) \setminus A$ such that $J_x \cap B_{2r_1}(z) = \emptyset$, we have
\begin{equation*}
  \harm_A^x (\cE_A(v)) \leq (1 - p)^K.
\end{equation*}
and the bound follows for $x = \infty$.
\end{proposition}

\begin{proof}
Items 1-4 in Definition~\ref{def:good-scales} identify $K$ concentric, disjoint, nonempty annuli in $G$ through which a random walk started outside of $\cV_{\cK}(B_{2r_1}(z))$ must pass, and for which a random walk started anywhere on the boundary has a disconnection probability of at least $p$.

In particular, item 2 states that $v$ and its immediate neighborhood are contained in the center of the concentric annuli, item 3 guarantees that a random walk trajectory disconnecting $\cV_{\cK}(\partial B_{3r_j/4}(z))$ from $\cV_{\cK}(\partial B_{5r_j/4}(z))$ must hit the cluster $A$, item 4 lower bounds the annulus intersection probability by $p$ and ensures that the visited plane data is contained in the annulus, and item 1 ensures that the random walk cannot skip scales and that the sets $\cV_{\cK}(\A_{r_j/3, 5r_j/3}(z))$ are disjoint for differing $j$.

As a random walk which first enters $A$ through $\cE_A(v)$ must avoid all $K$ disconnection events, one applies the Markov property of the random walk successively at each scale $j$ to conclude the final bound on the harmonic measure.
\end{proof}
To produce a diameter bound, we will need to show the existence of $(K, p)$-good scales uniformly over all $v \in A$. The next section shows how to convert such uniform good scales into a diameter bound.

\subsection{Harmonic measure bounds to diameter bounds}
We now use harmonic measure bounds to establish a diameter bound on the DLA cluster. Proposition~\ref{prop:helper-harm} operates at a fixed scale, but the next proposition takes as input harmonic measure bounds which tend to $0$ as the size of the subject cluster tends to $\infty$. Establishing these bounds at arbitrarily large scales, uniformly over clusters and vertices, is the main content of Sections~\ref{sec:pf-sublinear} and \ref{sec:pf-poly-sublinear}, and is done using the appropriate assumption in tandem with Proposition~\ref{prop:helper-harm}. We prove Proposition~\ref{prop:helper-diam} here since, as with Proposition~\ref{prop:helper-harm}, it does not inherently depend on either assumption.

Suppose that $G$ is a planar map with root $o = v_0$, and let $K_0 = \{o\}$, and $K_n$ the DLA cluster at $n$ steps, evolving according to Definition~\ref{def:dla}. $K_n$ has a tree structure, and we set $\textrm{height}(K_n) := \max_{v \in K_n} \dist(o,v;K_n)$. Define the (non-increasing) sequence
\begin{equation}\label{eq:harm-bound-asymp}
  a_n = a_n(G,o) := \sup_{\substack{A \ni o \text{ finite, connected} \\ |A| \geq n}}
  \ \sup_{v \in A} \harm_A^\infty(\cE_A(v)).
\end{equation}

Note that a priori the above quantity has nothing to do with DLA, and applies to all sufficiently large finite rooted connected subgraphs of $G$. The convergence of $a_n \to 0$, and in particular the rate of this convergence, is supplied by the Beurling estimates in Theorems~\ref{thm:beurling-weak} and \ref{thm:beurling-strong}; these latter results are established in subsequent sections, but here we show how to convert this harmonic measure control to internal diameter bounds on DLA.

\begin{proposition}\label{prop:helper-diam}
With $a_n$ as in \eqref{eq:harm-bound-asymp}, we have almost surely:
\begin{enumerate}
\item If $a_n \to 0$ then $\diam(K_n) = o(n)$.
\item If $a_n \leq C n^{-\theta}$ for some $C, \theta > 0$, then there exists a $\zeta > 0$ such that $\diam(K_n) = O(n^{1 - \zeta})$.
\end{enumerate}
\end{proposition}
\begin{proof}
By definition, for each $v \in K_{n - 1}$ the conditional probability that $e_n \in \cE_{K_{n - 1}}(v)$ given $K_{n - 1}$ is at most $a_n$, since $|K_{n - 1}| = n$. Let $s_n := \sum_{j = 1}^n a_j$. If $\textrm{height}(K_n) \geq k$, then there exist times $1 \leq t_1 < \dots < t_k \leq n$ such that $v_0,v_{t_1}, \dots, v_{t_k}$ form a simple path from the root. For any fixed increasing sequence $(t_1, \dots, t_k)$, successive conditioning on the growth history produces the bound $\prod_{j = 1}^k a_{t_j}$ on the event $\cap_{1 \leq j \leq k - 1} \{v_{t_j} \sim v_{t_{j + 1}} \} \bigcap \{v_0 \sim v_{t_1} \}$, where adjacency $\sim$ is in $K_n$. A union bound over all such sequences gives $\P[\, \textrm{height}(K_n) \geq k] \leq \sum_{1 \leq t_1 < \cdots < t_k \leq n} \prod_{j = 1}^k a_{t_j} \leq s_n^k/k! \leq (e s_n/k)^k$. The second inequality follows by expanding $s_n^k$, and the last follows from $k! \geq (k/e)^k$.

First, suppose that $a_n \to 0$. Then $s_n/n \to 0$ and so $s_n = o(n)$. Fix $c > 0$ and take $k = \lceil cn \rceil$. For all sufficiently large $n$ one has $e s_n/\lceil c n \rceil \leq 1/2$, and so $\P[\textrm{height}(K_n) \geq \lceil cn \rceil] \leq 2^{-\lceil cn \rceil}$. Since $\sum_{n \in \N} \P[\textrm{height}(K_n) \geq \lceil cn \rceil] < \infty$, by Borel-Cantelli we have $\textrm{height}(K_n) \leq O(cn)$ a.s., with an implied constant independent of $c$. By intersecting these probability-one events over all positive rational $c$, we conclude $\textrm{height}(K_n) = o(n)$ a.s.

Now suppose that $a_n \leq C n^{-\theta}$ and let $\theta_0 := \min \{ \theta, 1/2\}$ (or $\theta_0 := \min \{ \theta, c\}$ for any fixed $c \in (0, 1)$). We have $s_n \leq C'n^{1 - \theta_0}$ (with possibly $C' > C$). Choose $\zeta_0 \in (0, \theta_0)$ and set $k_n := \lceil n^{1 - \zeta_0} \rceil$. We then have $\P[\textrm{height}(K_n) \geq k_n] \leq (C''n^{\zeta_0 - \theta_0})^{n^{1 - \zeta_0}}$ and Borel-Cantelli implies that a.s.\ $\textrm{height}(K_n) = O(n^{1 - \zeta_0})$.

Since $\diam(K_n) \leq 2 \cdot \textrm{height}(K_n)$, the claimed diameter bounds follow in both cases.
\end{proof}

Note that $K_n$, viewed not as a tree but with all edges from $G$ between vertices of $K_n$ retained, trivially has diameter at most that of the tree, so conforms to the asymptotics of Proposition~\ref{prop:helper-diam} too. 

To invoke Proposition~\ref{prop:helper-diam}, we must show existence of good scales for the DLA cluster $K_n$ for all sufficiently large $n$, with the number of good scales growing in $n$, after which we can use Proposition~\ref{prop:helper-harm}.

\subsection{Change of initial vertex}
We now briefly describe the procedure for changing the starting vertex for DLA growth. For $b \in \cV(G)$, define $a_n(G,b)$ by replacing $o$ with $b$ in~\eqref{eq:harm-bound-asymp}. Fix a finite path $P$ from $o$ to $b$ and let $T_S := \inf \{j \geq 0:X_j \in S\}$ for the random walk $X_n$. Set
\begin{equation*}
  q_P := \min_{u \in P} \P_u(T_{\{b\}} < \infty) > 0.
\end{equation*}
For finite connected $A \ni b$, put $B = A \cup P$. A walk first entering $A$ at $v \in A$ must first enter $B$ at $v$ or in $P \setminus A$. Also, for $x \notin B$ the Markov property gives $\P_x(T_A < \infty) \geq q_P \P_x(T_B < \infty)$. Conditioning on hitting the respective sets and then sending $x$ to infinity gives
\begin{align*}
  \harm_A^\infty(\cE_A(v)) \leq q_P^{-1}(\harm_B^\infty(\cE_B(v))
       + \sum_{u \in P \setminus A} \harm_B^\infty(\cE_B(u))) \leq \frac{|P|}{q_P}a_{|A|}(G,o),
\end{align*}
since $B$ is rooted and connected, $|B| \geq |A|$ and $1 + |P \setminus A| \leq |P|$. Thus
\begin{equation}\label{eq:root-comparison}
  a_n(G,b) \leq \frac{|P|}{q_P}a_n(G,o).
\end{equation}
For a recurrent random walk we have $q_P = 1$, and the same argument applies to the conductance weighted walk. Proposition~\ref{prop:helper-diam} therefore gives the same growth bounds for DLA started from $b$, with the same polynomial exponent where applicable. However, the implied constants and eventual thresholds may depend on the initial vertex, and so change of initial growth point for an RPM is a little more subtle. For this reason, we largely omit the discussion of change of initial growth point in the RPM case, and focus on the rooted growth.

\section{Proof of Theorem~\ref{thm:sublinear} under Assumption~\ref{assumption:rw-conv}}\label{sec:pf-sublinear}
We now use Assumption~\ref{assumption:rw-conv} to show the existence of good scales, and invoke Propositions~\ref{prop:helper-harm} and \ref{prop:helper-diam}. First, we define the relevant plane data $\cK$ in this setting. For an edge $e$, let $\Phi(e)^\circ$ denote the embedded edge with endpoints removed. For $x \in \cV(G)$, let
\begin{equation*}
  J_x := \{ \Phi(x)\} \cup \bigcup_{\substack{e \in \cE(G): \, x \in e}} \Phi(e)^\circ.
\end{equation*}
Each $J_x$ is bounded and connected, and the family $\cK := \{J_x: x \in \cV(G) \}$ is locally finite by properness of the embedding (see Definition~\ref{def:embedding}). Moreover, for distinct vertices $x,y$,
\begin{equation*}
  J_x \cap J_y \neq \emptyset \quad \Longleftrightarrow \quad x \sim y.
\end{equation*}
As before, for $D \subset \C$ we write
\begin{equation*}
  \cV_\Phi(D) := \{x \in \cV(G):J_x \cap D \neq \emptyset \}, \quad 
  D(A) := \diam \{ \Phi(x): x \in A\}.
\end{equation*}
Note that generally $D(A) \neq \diam(U_A^{\cK})$.

\subsection{Robust Brownian intersections}\label{subsec:robust-brownian-intersections}
Here, we prove that with positive probability a Brownian motion intersects itself in an annulus in such a way that an appropriate approximating curve, which we will take to be a random walk trajectory on the embedded graph, also intersects itself. This allows us to infer the relevant annular intersection probability of the random walk from the Brownian motion convergence. The claims of this section are elementary and are therefore stated without proof.

Let $\Gamma$ be the polygonal path depicted in the below figure. The first vertical segment and the last horizontal segment cross transversely at $(1,1/10)$, and the loop winds once around the origin while staying in a small annulus.

\begin{figure}[ht]
\centering

\begin{tikzpicture}[
    scale=3.2,
    line cap=round,
    line join=round,
    route/.style={line width=1.1pt,-{Stealth[length=2.1mm,width=1.5mm]}},
    vertex/.style={circle,fill=black,inner sep=1.05pt},
    label/.style={font=\scriptsize,fill=white,fill opacity=.9,
                  text opacity=1,inner sep=1.2pt},
    outer annulus/.style={red,densely dotted,line width=.75pt},
    inner annulus/.style={blue,densely dotted,line width=.75pt},]

    \draw[outer annulus] (0,0) circle[radius={3/4}];
    \draw[outer annulus] (0,0) circle[radius={5/4}];
    \draw[inner annulus] (0,0) circle[radius={4/5}];
    \draw[inner annulus] (0,0) circle[radius={6/5}];
    
    \coordinate (p0)  at (1,0);
    \coordinate (p1)  at (1,{1/4});
    \coordinate (p2)  at ({3/4},{3/4});
    \coordinate (p3)  at (0,1);
    \coordinate (p4)  at ({-3/4},{3/4});
    \coordinate (p5)  at (-1,0);
    \coordinate (p6)  at ({-3/4},{-3/4});
    \coordinate (p7)  at (0,-1);
    \coordinate (p8)  at ({3/4},{-3/4});
    \coordinate (p9)  at (1,{-1/4});
    \coordinate (p10) at ({9/10},{1/10});
    \coordinate (p11) at ({11/10},{1/10});
    \coordinate (q)   at (1,{1/10});
    
    \draw[route]
      (p0)--(p1)--(p2)--(p3)--(p4)--(p5)--
      (p6)--(p7)--(p8)--(p9)--(p10)--(p11);
    
    \foreach \p in {p0,p1,p2,p3,p4,p5,p6,p7,p8,p9,p10,p11}
      \node[vertex] at (\p) {};
    
    \draw[fill=white,line width=.45pt] (q) circle[radius=.014];
    
    \node[label,anchor=north east,xshift=-2pt,yshift=-2pt] at (p0)
      {$\left(1,0 \right)$};
    \node[label,anchor=south west,xshift=2pt,yshift=2pt] at (p1)
      {$\left(1, \frac{1}{4} \right)$};
    \node[label,anchor=south west,xshift=2pt,yshift=2pt] at (p2)
      {$\left(\frac{3}{4}, \frac{3}{4} \right)$};
    \node[label,anchor=south,yshift=2pt] at (p3)
      {$\left(0,1 \right)$};
    \node[label,anchor=south east,xshift=-2pt,yshift=2pt] at (p4)
      {$\left(-\frac{3}{4}, \frac{3}{4} \right)$};
    \node[label,anchor=east,xshift=-2pt] at (p5)
      {$\left(-1,0 \right)$};
    \node[label,anchor=north east,xshift=-2pt,yshift=-2pt] at (p6)
      {$\left(-\frac{3}{4},-\frac{3}{4} \right)$};
    \node[label,anchor=north,yshift=-2pt] at (p7)
      {$\left(0,-1 \right)$};
    \node[label,anchor=north west,xshift=2pt,yshift=-2pt] at (p8)
      {$\left(\frac{3}{4},-\frac{3}{4} \right)$};
    \node[label,anchor=north west,xshift=2pt,yshift=-2pt] at (p9)
      {$\left(1,-\frac{1}{4} \right)$};
    \node[label,anchor=south east,xshift=-2pt,yshift=3pt] at (p10)
      {$\left(\frac{9}{10}, \frac{1}{10} \right)$};
    \node[label,anchor=west,xshift=3pt] at (p11)
      {$\left(\frac{11}{10}, \frac{1}{10} \right)$};
\end{tikzpicture}

\vspace{0.6em}

\caption{The polygonal path $\Gamma$. The open circle at $(1, \frac{1}{10})$ indicates the initial intersection of the first vertical segment with the final horizontal segment. The dotted blue circles indicate the annulus $\A_{4/5, 6/5}(0)$, and the dotted red circles indicate the annulus $\A_{3/4, 5/4}(0)$.}
\label{fig:polygonal-path}
\end{figure}
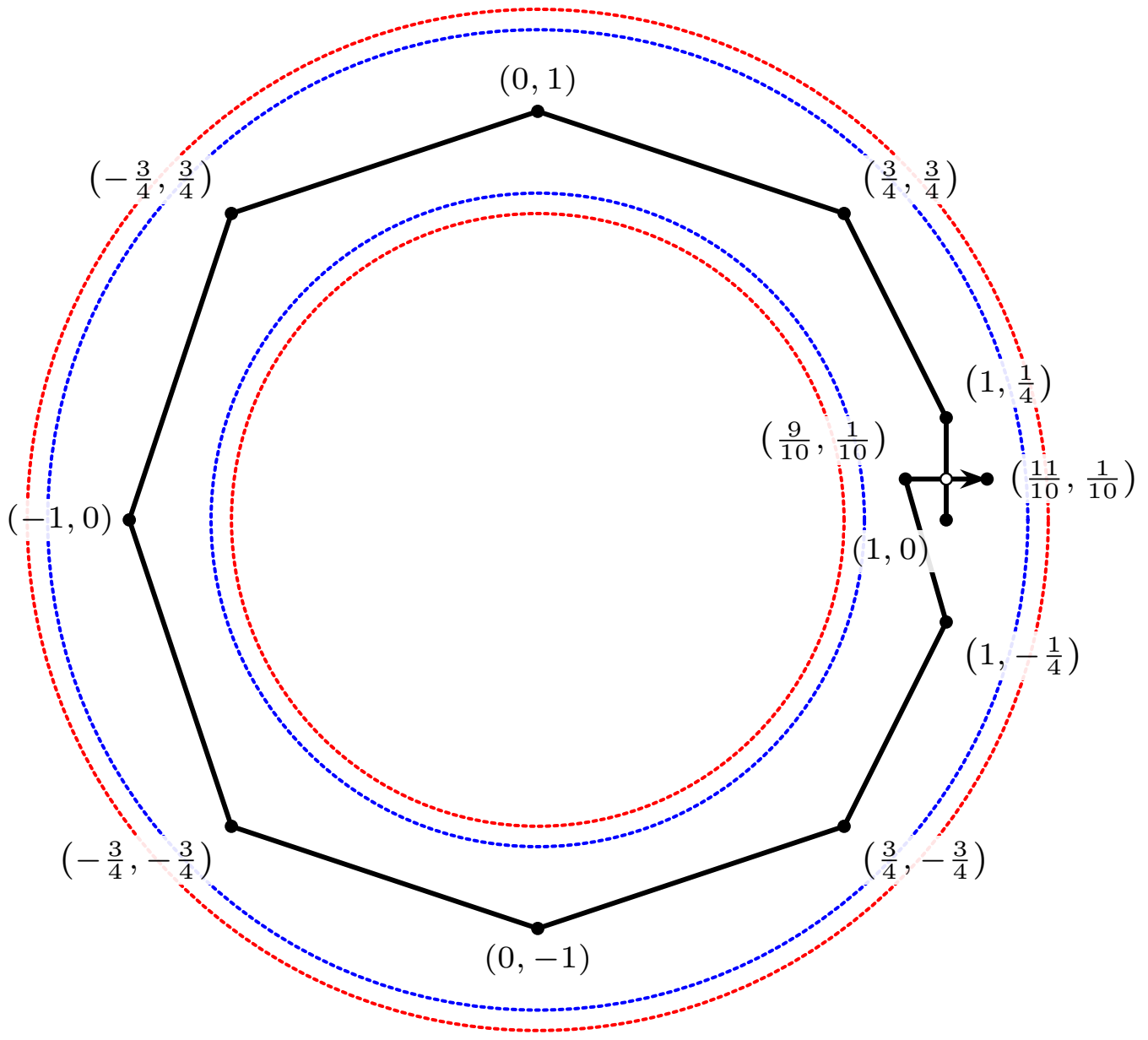

The next claim is obvious.

\begin{claim}\label{claim:robust-intersection}
There is $\eta_* > 0$ such that every finite curve $\gamma$ within $\eta_*$ of $\Gamma$ in the metric $d_F$ stays in $\A_{3/4, 5/4}(0)$ and contains a closed curve with winding number one about $\overline{B_{3/4}(0)}$.
\end{claim}

We call an $\epsilon$-neighborhood of $\Gamma$ with $\epsilon < \eta_*$ a \textit{robust neighborhood}. The next claim is immediate by planarity and properties of the embedding.

\begin{claim}\label{claim:cycles}
If a finite piece of an embedded walk contains a closed curve with nonzero winding number about $z$, then the graph walk contains an embedded graph cycle surrounding $z$.
\end{claim}

The final claim of this section is a basic property of Brownian motion (see e.g. \cite[Exercise 1.8]{morters2010brownian}).

\begin{claim}\label{claim:brownian-robust-intersection}
There exist $p_0 > 0$ and $\eta \in (0, \eta_*/10)$ such that for $y = \rho e^{i \theta}$ with $1 - \eta \leq \rho \leq 1 + \eta$, standard Brownian motion started from $y$ has probability at least $p_0$ to possess an initial segment in a robust neighborhood of $\Gamma_y$ (where $\Gamma_y$ is the rotation by $\theta$ of the template $\Gamma$ whose initial point $(1,0)$ has been replaced by $(\rho,0)$).
\end{claim}

We allow $\rho \neq 1$ because, in our later applications to embedded random walks, the starting vertex will generally lie close to a prescribed Euclidean circle rather than exactly on it. We call the event that Brownian motion, and by abuse of terminology the random walk too, is in a robust neighborhood the \textit{robust intersection event}.

\subsection{Diameter bounds for embedded edges}
As mentioned in Section~\ref{subsec:annular-intersections}, to ensure that the embedded graph conforms to the non-normalized scales, we need to bound the diameter of the elements in the planar data. In this subsection, we prove that the maximum diameter of the edges (of a graph satisfying Assumption~\ref{assumption:rw-conv}) intersecting a Euclidean ball is sublinear in the radius of the ball.

For $C, R > 0$ define
\begin{equation*}
  m_C(R) := \frac{1}{R} \max \{ \diam(\Phi(e)):e \in \cE(G),\, \Phi(e) \cap \overline{B_{CR}(0)} \neq \emptyset \}.
\end{equation*}
We have $m_C(R) < \infty$ by properness of the embedding.

\begin{proposition}\label{prop:sublinear-mesh}
Under Assumption~\ref{assumption:rw-conv} for every fixed $C < \infty$ we have
\begin{equation*}
  m_C(R) \to 0 \text{ as } R \to \infty.
\end{equation*}
\end{proposition}
\begin{proof}
Fix $h \in (0,1)$. We prove that $m_C(R) \leq 5h/2$ for all sufficiently large $R$. Choose finitely many points $z_1, \dots, z_N \in \overline{B_C(0)}$ such that
\begin{equation}
    \overline{B_C(0)}
    \subset \bigcup_{i = 1}^N B_{h/4}(z_i).
   \label{eq:finite-cover}
\end{equation}
For each $i$, choose a polygonal curve which starts at $0$, travels to $z_i + (h,0)$, traces a translated and scaled copy of $\Gamma$, and then continues to $\partial B_L(0)$. We choose $L > C + 10$ so that the polygonal curves and their relevant neighborhoods are contained in $B_L(0)$.

Let $d_F$ be the metric used in Assumption~\ref{assumption:rw-conv} and $\W$ be the law of Brownian motion started at $0$ and stopped upon exiting $B_L(0)$. Additionally, let $\mu_R$ be the law of the rescaled embedded random walk started from the root $o$ and stopped upon leaving $B_L(0)$. Assumption~\ref{assumption:rw-conv} gives
\begin{equation*}
    d_{\textrm{Pr}}(\W, \mu_R) \to 0 \quad \text{ as } R \to \infty.
\end{equation*}
For each $i$ we let $\gamma_i$ be the corresponding polygonal route. By the Claim~\ref{claim:robust-intersection} there is $\eta_i > 0$ such that every curve $\gamma$ satisfying $d_{F}(\gamma, \gamma_i) < 2 \cdot \eta_i$ contains a curve in $\A_{3h/4, 5h/4}(z_i)$ separating $\overline{B_{3h/4}(z_i)}$ from $\infty$. Define the events
$$
    \cA_i := \{ \gamma:d_{F}(\gamma, \gamma_i) < \eta_i\},
    \qquad
    \cB_i := \{ \gamma:d_{F}(\gamma, \gamma_i) < 2 \cdot \eta_i\}.
$$
Claim~\ref{claim:brownian-robust-intersection} gives $q_i := \W(\cA_i) > 0$, and since $N < \infty$, the number $\epsilon_{C,h} := \min_{1 \leq i \leq N}\{ \eta_i, q_i/2\}$ is positive. By Assumption~\ref{assumption:rw-conv} there is a finite $R_*(C,h)$ such that
\begin{equation*}
     d_{\textrm{Pr}}(\W, \mu_R) < \epsilon_{C,h} \qquad \text{for every } R \geq R_*(C,h).
\end{equation*}
If $d_{\textrm{Pr}}(\W, \mu_R) < \epsilon_{C,h}$ then by definition of the L\'evy-Prokhorov metric we have for every Borel set $\cA$ (in the topology induced by $d_F$) that
$$
    \W(\cA) \leq \mu_R(\cA^{\epsilon_{C,h}}) + \epsilon_{C,h},
$$
where $\cA^\epsilon$ denotes the open $\epsilon$-neighborhood of $\cA$ under $d_F$. Taking $\cA = \cA_i$ gives $\mu_R(\cA_i^{\epsilon_{C,h}}) \geq q_i - \epsilon_{C,h} \geq q_i/2 > 0.$ Moreover, since $\epsilon_{C,h} \leq \eta_i$ we have $\cA_i^{\epsilon_{C,h}} \subset \cB_i$ and so
\begin{equation}
    \mu_R(\cB_i) \geq \frac{q_i}{2} > 0
    \qquad
    \text{for every $i$ and every $R \geq R_*(C,h)$}.
   \label{eq:positive-discrete-lasso}
\end{equation}
Fix $R \geq R_*(C,h)$. Since the graph and embedding are deterministic, and the randomness in \eqref{eq:positive-discrete-lasso} is only over the random walk, positive probability implies that an appropriate trajectory exists. Choose one such trajectory for each $i$. By Claims~\ref{claim:robust-intersection} and~\ref{claim:cycles} each trajectory contains a graph cycle $\mathcal{C}_i^R$ such that
\begin{equation*}
    R^{-1} \Phi(\mathcal{C}_i^R) \subset \A_{3h/4, 5h/4}(z_i) \text{ and } R^{-1} \Phi(\mathcal{C}_i^R) \text{ surrounds } \overline{B_{3h/4}(z_i)}.
\end{equation*}
Since the above probabilistic argument is only existential, we can deduce the existence of such cycles for all $i$ simultaneously. Let $e$ be an edge such that $\Phi(e) \cap \overline{B_{CR}(0)} \neq \emptyset$. By \eqref{eq:finite-cover}, there is an $i$ such that
$$
    R^{-1} \Phi(e) \cap B_{h/4}(z_i) \neq \emptyset.
$$
The edge $e$ cannot belong to $\mathcal{C}_i^R$ since the cycle is disjoint from the disk, and by planarity the interior of $e$ cannot cross the cycle. So $R^{-1} \Phi(e)$ is contained in $B_{5h/4}(z_i)$, and $\diam(R^{-1} \Phi(e)) \leq \frac{5}{2}h.$ This bound holds for all edges, and so for each fixed $h > 0$ we have $\limsup_{R \to \infty} m_C(R) \leq \frac{5}{2}h$ and letting $h \downarrow 0$ completes the proof.
\end{proof}

It follows immediately that for every fixed $C < \infty$ one has $\frac{1}{R} \max \{ \diam(J_x): J_x \cap \overline{B_{CR}(0)} \neq \emptyset \} \to 0$ as $R \to \infty$.

\subsection{Completing Theorem~\ref{thm:sublinear}}
Now, we verify the existence of good scales and complete the proof of Theorems~\ref{thm:beurling-weak} and \ref{thm:sublinear}. The next Lemma establishes control of annular intersection probabilities at a single normalized scale $\rho \in (0, 1)$.

\begin{lemma}\label{lem:single-scale}
Under Assumption~\ref{assumption:rw-conv} there exist deterministic $\eta, p \in (0, 1)$ such that for every $C < \infty$ and every fixed normalized scale $\rho \in (0,1)$, there is $R_*(C, \rho) < \infty$ such that the following holds for all $R \geq R_*(C, \rho)$. Let $z \in B_{CR}(0)$, set $r := \rho R$ and let $x \in \cV(G)$ satisfy
\begin{equation}\label{eq:admissable-starting-point}
    1 - \eta \leq \frac{|\Phi(x) - z|}{r} \leq 1 + \eta.
\end{equation}
 A simple random walk $X_n$ started from $x$ and run until its first exit time $\tau$ from $\cV_{\cK}(\A_{r/2, 3r/2}(z))$ has probability at least $p$ to disconnect $\cV_{\cK}(\partial B_{3r/4}(z))$ from $\cV_{\cK}(\partial B_{5r/4}(z))$ with $\cup_{m = 0}^{\tau} \, J_{X_m} \subset \A_{r/3, 5r/3}(z)$. The probability $p$ is independent of $C$, $\rho$, $R$, $z$, and $x$.
\end{lemma}
\begin{proof}
Fix $C < \infty$ and $\rho \in (0,1)$, and set $r = \rho \cdot R$. Choose $L > C + 3$, and let $\eta$ be as in Claim~\ref{claim:brownian-robust-intersection}, decreasing $\eta$ if necessary such that $\eta < 1$. By \eqref{eq:admissable-starting-point} every admissible starting vertex $x$ satisfies $|\Phi(x)| \leq (C + 2)R$. After translating by $-z$, rescaling by $r^{-1}$ and rotating so that the starting point lies near $(1,0)$, Claim~\ref{claim:brownian-robust-intersection} gives a $p_0 > 0$, independent of $C, \rho,z$, and $x$ such that Brownian motion remains in a robust neighborhood of $\Gamma$. In the original coordinates rescaled by $R^{-1}$, this curve remains in $\A_{3 \rho/4, 5 \rho/4}(z/R)$ and contains a closed curve in $\A_{4 \rho/5, 6 \rho/5}(z/R)$ which separates $\overline{B_{3 \rho/4}(z/R)}$ from $\infty$. By an identical argument to that of Proposition~\ref{prop:sublinear-mesh} there exists an $R_1(C, \rho) < \infty$ such that when $R \geq R_1(C, \rho)$ the rescaled embedded random walk intersects robustly in the relevant annulus uniformly over all admissible $z$ and $x$ with probability at least $p_0/2$.

By Proposition~\ref{prop:sublinear-mesh} there is an $R_2(C, \rho) < \infty$ such that for $R \geq R_2(C, \rho)$, every set $J_u$ intersecting the ambient region containing the robust neighborhood of (the rotation of) $\Gamma$ has diameter smaller than a sufficiently small multiple of $r$. On the robust intersection event, by Claim~\ref{claim:cycles} the trace of the random walk contains a cycle in $\A_{4r/5, 6r/5}(z)$ disconnecting $\overline{B_{3r/4}(z)}$ from $\infty$.

The containment claims of the Lemma follow analogously to those of Proposition~\ref{prop:sublinear-mesh}, and choosing $R_*(C, \rho) = \max (R_1(C, \rho), R_2(C, \rho))$ completes the proof.
\end{proof}
To observe that the above $p$ is independent of $\rho$ we use the scale invariance of Brownian motion modulo parametrization and the fact that the Brownian motion is killed upon leaving the annulus, which is itself appropriately rescaled.

Now, we show the existence of many good scales as the Euclidean diameter of the cluster diverges. The next proposition essentially uses Proposition~\ref{prop:sublinear-mesh} to show the existence of an arbitrarily large number of disjoint scales as $R \to \infty$ and then applies Lemma~\ref{lem:single-scale} at each individual scale. Uniformity over $A$ and $v$ follows from the fact that the sets are rooted. 

\begin{proposition}\label{prop:many-good-scales}
Under Assumption~\ref{assumption:rw-conv} there is a $p > 0$ such that for every $K \in \N$ there exists an $m_K < \infty$ with the following property: if $A \ni o$ is finite and connected, $|A| \geq m_K$ and $v \in A$, then there are $(K,p)$-good scales for $(A,v)$.
\end{proposition}
\begin{proof}
Let
\begin{equation*}
  \rho_1 := 1/100, \qquad \rho_j := \rho_1 \cdot 64^{-(j - 1)}, \quad 1 \leq j \leq K.
\end{equation*}
Properness of the embedding therefore implies that $\inf_{A \ni o: |A| > k}D(A) \to \infty$ as $k \to \infty$ over rooted connected sets; thus, any property which holds for $R$ large enough also holds for large enough $k$.

Now, suppose that $K \in \mathbb{N}$ and we choose $\rho_1$ and $\rho_K$ as above. Since the maximum diameter of all $J_x$ intersecting $B_{2R}$ is $o(R)$, by choosing $C = 2$ in Proposition~\ref{prop:sublinear-mesh} there exists an $R(K)$ (where $\eta$ and $p$ are as in Lemma~\ref{lem:single-scale}) such that if $R \geq R(K)$ then the maximum diameter of all $J_x$ intersecting $B_{2R}$ is less than $\rho_K R \min \{1/100, \eta \}$.

For $A$ a rooted, finite, connected set, let $R = D(A)$ and fix $z = \Phi(v)$. Since $A$ is rooted and $\Phi(o) = 0$, we have $\Phi(A) \subset \overline{B_R(0)}$. For every $x \in N(v) \setminus \{v\}$, $J_x \cap J_v \neq \emptyset$, so for $R \geq R(K)$, each $J_x$ intersects $B_{2R}(0)$ and has diameter less than $\rho_K R/100$. If $D(A) \geq R(K)$ then $\bigcup_{x \in N(v)}J_x$ is contained in $B_{\rho_K R/20}(z)$, and every $J_x$ intersecting $B_{2 \rho_1 R}(z)$ has diameter less than $\rho_K R/100$ by the preceding paragraph, since $B_{2 \rho_1 R}(z) \subset B_{2R}(0)$.

There exists $u \in A$ with $|\Phi(u) - z| \geq R/2$, and an embedded graph path in $A$ from $v$ to $u$ crosses both circles $\partial B_{3 \rho_j R/4}(z)$ and $\partial B_{5 \rho_j R/4}(z)$ for all $1 \leq j \leq K$. Applying Lemma~\ref{lem:single-scale} with $C = 2$ to each $\rho_j$ shows the existence of $(K, p)$-good scales, and picking $m_K$ such that $\inf_{A \ni o: |A| \geq m}D(A) > R(K)$ for all $m \geq m_K$ completes the proof.
\end{proof}

Note that the order of quantifiers to achieve uniformity over $A$ large and $v$ is somewhat delicate. For the fixed threshold $R(K)$, an $A$ with $D(A) > R(K)$ may fail to satisfy $\max \{ \diam(J_x):x \in A\} \leq \rho_K R(K) /100$, since there may exist some $v \in A$ such that $J_v \cap B_{2R(K)}(0) = \emptyset$. However, $\max \{ \diam(J_x):x \in A\} \leq \rho_K D(A) /100$ does hold, since the numerator in the definition $m_C(R)$ adjusts the size of the ball with which edges intersect, and so the convergence $m_C(R) \to 0$ controls the relative mesh size at all sufficiently large scales simultaneously.

Now, we can complete Theorems~\ref{thm:beurling-weak} and \ref{thm:sublinear}.

\begin{proof}[Proof of Theorems~\ref{thm:beurling-weak} and \ref{thm:sublinear}]
Let $p$ be as in Proposition~\ref{prop:many-good-scales}. Given $\delta > 0$, choose $K$ with $(1 - p)^K < \delta$. For every sufficiently large rooted connected $A$ and every $v \in A$ there are $(K,p)$-good scales for $(A,v)$. By Proposition~\ref{prop:helper-harm} we have
$$
  \harm_A^\infty(\cE_A(v)) \leq (1 - p)^K < \delta.
$$
Sending $K \to \infty$ as $\delta \to 0$ produces Theorem~\ref{thm:beurling-weak}, and combining with Proposition~\ref{prop:helper-diam} completes Theorem~\ref{thm:sublinear} for DLA started from $o$. Equation~\eqref{eq:root-comparison} gives the extension to any other initial vertex.
\end{proof}

\section{Proof of Theorem~\ref{thm:poly-sublinear} under Assumption~\ref{assumption:mesoscopic}}\label{sec:pf-poly-sublinear}
As in Section~\ref{sec:pf-sublinear}, we use Assumption~\ref{assumption:mesoscopic} to show the existence of good scales, and invoke Proposition~\ref{prop:helper-diam}. The stronger control of Assumption~\ref{assumption:mesoscopic} as compared to Assumption~\ref{assumption:rw-conv} will allow us to improve the number of good scales.

We begin with two trivial observations. The first follows from the triangle inequality, and the second follows from the reverse triangle inequality and the connectedness of the union of cells corresponding to a connected subgraph of $G$.

\begin{lemma}\label{lem:neighbor-local}
Suppose that for every $x \in N(v)$ the cell $H_x$ has diameter at most $d$, and choose $z_v \in H_v$. Then
\begin{equation*}
  \bigcup_{x \in N(v)}H_x \subset \overline{B_{2d}(z_v)}.
\end{equation*}
\end{lemma}

\begin{lemma}\label{lem:cluster-crossing}
Let $A$ be finite and connected, with $D := \diam(U_A^{\cH})$, and suppose $a > 0$ and $\dist(z,U_A^{\cH}) \leq a < b \leq D/2$. Then $A$ contains a graph path joining $\cV_{\cH}(\partial B_a(z))$ to $\cV_{\cH}(\partial B_b(z))$.
\end{lemma}

The next proposition establishes good scales for every sufficiently large non-normalized scale, where the number of good scales depends logarithmically on the Euclidean diameter of the cluster. Subsequently, item 4 of Assumption~\ref{assumption:mesoscopic} will be used to convert the lower bound on the number of good scales from a term depending on the Euclidean diameter of the cells corresponding to the cluster to a term depending only on the cardinality of the cluster.

\begin{proposition}\label{prop:mesoscopic-good-scales}
Under Assumption~\ref{assumption:mesoscopic}, there exists a constant $c_0 > 0$ and a finite threshold $R_1 \geq R_0$ such that the following holds: if $A \ni o$ is finite and connected, $R := \diam(U_A^{\cH}) \geq  R_1$, and $v \in A$, then there are $(K,p)$-good scales for $(A,v)$ with
\begin{equation*}
  K \geq c_0 \log R.
\end{equation*}
\end{proposition}

\begin{proof}
Fix $C_0 \geq 200$. By Assumption~\ref{assumption:mesoscopic}, there exists an $R_1 \geq R_0$ sufficiently large that for every $R \geq R_1$ we have
\begin{equation}\label{eq:mesoscopic-inequalities}
  e(R) + d(R) \leq m(R), \qquad C_0m(R) \leq R/200.
\end{equation}
Set $R := \diam(U_A^{\cH}) \geq R_1$, fix $v \in A$, and choose $z_v \in H_v$. Since $A$ is rooted and $0 \in H_o$, we have $0 \in U_A^{\cH}$ and hence $U_A^{\cH} \subset \overline{B_R(0)}$. Next, choose $z \in \mathcal{Z}_R$ with $|z - z_v| \leq e(R)$. Set
\begin{equation*}
  s_j := \frac{R}{200}32^{-(j - 1)}, \qquad j \geq 1,
\end{equation*}
and let $K$ be the largest integer with $s_K \geq C_0m(R)$. For each $j \leq K$, choose $r_j \in \mathcal{Q}_R$ with $|r_j - s_j| \leq e(R)$. Equation~\eqref{eq:mesoscopic-inequalities} gives $e(R) \leq m(R) \leq s_j/C_0$, and so we have $r_j \geq (1 - C_0^{-1})s_j$ and $r_{j + 1} \leq (32^{-1} + C_0^{-1})s_j$, so $r_{j + 1} \leq r_j/8$. By our choice of $K$, we have $K > \log_{32}(R/(200C_0m(R)))$. Since $m(R) \leq R^{1 - \sigma}$, by possibly increasing $R_1$ again, we have
\begin{equation*}
  K \geq c_0 \log R.
\end{equation*}

Now, we can verify the existence of $(K,p)$-good scales for $(A,v)$. Since $z \in B_{2R}(0)$ and $r_1 \leq R/100$, all cells intersecting $B_{2r_1}(z)$ also intersect $B_{4R}(0)$, so by Assumption~\ref{assumption:mesoscopic} the cell diameters are bounded by $d(R)$. Lemma~\ref{lem:neighbor-local} then gives $\cup_{x \in N(v)} H_x \subset \overline{B_{2d(R) + e(R)}(z)}$. Now, $r_K \geq s_K - e(R) \geq (C_0 - 1)m(R)$, so $d(R) \leq r_K/100$ and $2d(R) + e(R) \leq r_K/20$, thus
\begin{equation*}
  \bigcup_{x \in N(v)}H_x \subset B_{r_K/20}(z).
\end{equation*}
This proves the first two items in Definition~\ref{def:good-scales}. Since $z_v \in U_A^{\cH}$ we have $\dist(z,U_A^{\cH}) \leq e(R) < 3r_K/4$. Also, $5r_j/4 \leq 5r_1/4 < R/2$ after increasing $R_1$ again if necessary. Applying Lemma~\ref{lem:cluster-crossing} with $a = 3r_j/4$ and $b = 5r_j/4$ gives the third item in Definition~\ref{def:good-scales} for every $R$. Finally, $z \in \mathcal{Z}_R$ and $r_j \in \mathcal{Q}_R$, and so the third item in Assumption~\ref{assumption:mesoscopic} gives the fourth item in Definition~\ref{def:good-scales}. This concludes the proof.
\end{proof}

We can now conclude Theorems~\ref{thm:beurling-strong} and \ref{thm:poly-sublinear}.

\begin{proof}[Proof of Theorems~\ref{thm:beurling-strong} and \ref{thm:poly-sublinear}]
Let $A \ni o$ be finite and connected with $\diam(U_A^{\cH}) \geq R_1$. Propositions~\ref{prop:helper-harm} and \ref{prop:mesoscopic-good-scales} give~
\begin{equation*}
  \harm_A^\infty(\cE_A(v))
  \leq (1 - p)^{c_0 \log \diam(U_A^{\cH})}
  = \diam(U_A^{\cH})^{-\xi},
  \qquad \xi := -c_0 \log(1 - p) > 0.
\end{equation*}
By item four of Assumption~\ref{assumption:mesoscopic}, for $A$ sufficiently large we have
\begin{equation*}
  \diam(U_A^{\cH})^{-\xi}
  \leq c_{\vol}^{-\xi}|A|^{-\nu \xi}.
\end{equation*}
Choosing any $\theta \in (0, \nu \xi)$ gives the strong Beurling estimate for sufficiently large $A$. Theorem~\ref{thm:beurling-strong} and item two of Proposition~\ref{prop:helper-diam} give the polynomially sublinear growth bound for DLA started from $o$. Equation~\eqref{eq:root-comparison} transfers this bound to any other initial vertex with the same exponent $\zeta$.
\end{proof}

\section{Verifying Assumption~\ref{assumption:mesoscopic} for the mated CRT map}\label{sec:mated-verification}
With the preliminaries from Section~\ref{sec:prelims} we verify Assumption~\ref{assumption:mesoscopic} for the mated CRT map, and conclude Theorem~\ref{thm:dla-mated-crt}. As the assumption pertains to a deterministic geometry, we mean that the mated CRT map satisfies the assumption a.s.\ and with possibly random values for constants such as $R_0$.

We will need to convert the properties stated in Proposition~\ref{prop:mated-crt-mesoscopic}, which pertain to a single scale, to the multiscale requirements of Assumption~\ref{assumption:mesoscopic}. Essentially, the argument boils down to scaling properties of the mated CRT map and an elementary Borel-Cantelli argument.

\begin{proposition}\label{prop:mated-mesoscopic}
Fix $\gamma \in (0,2)$. The $\gamma$-mated CRT map, rooted at $o = 0$, a.s.\ satisfies Assumption~\ref{assumption:mesoscopic}, with deterministic, possibly $\gamma$-dependent, $p$ and $\sigma$.
\end{proposition}
\begin{proof}
In what follows we work on a full probability event over the mated CRT map, so omit the phrase ``almost surely'' as regards the randomness of the environment.

Let $p, \alpha, \beta, \chi,q, \rho$ be as in Proposition~\ref{prop:mated-crt-mesoscopic}, decreasing $\beta$ if necessary such that $\beta < a_\gamma$. For $k \in \N$, set $b_k := 2^k$, $m_k := R_{b_k} b_k^{-\beta}$, $e_k := R_{b_k} b_k^{-\chi}$, and $d_k := R_{b_k} b_k^{-q}$.

By quantum cone scaling (Proposition~\ref{prop:quantum-cone-scaling}) the conclusion of Proposition~\ref{prop:mated-crt-mesoscopic} at $\epsilon = 2^{-k}$ transfers to the unit-mesh map at the (random) non-normalized scale $R_{b_k}$. By Borel-Cantelli, simultaneously for all sufficiently large $k$ the conclusions of Proposition~\ref{prop:mated-crt-mesoscopic} hold, i.e. every cell intersecting $B_{\rho R_{b_k}}(0)$ has diameter at most $d_k$ and the annular random walk conclusions of Proposition~\ref{prop:mated-crt-mesoscopic} are true.

Choose $\lambda \in (1/a_\gamma,1/(a_\gamma - \beta))$ and set $\sigma_* := 1 - \lambda(a_\gamma - \beta) > 0$. For $D > 1$, let $k(D) := \lceil \lambda \log_2D \rceil$; we will treat $D$ as the non-normalized scale parameter. Proposition~\ref{prop:radius-asymp} gives $R_{b_{k(D)}}/D \to \infty$, $m_{k(D)} = D^{1 - \sigma_* + o_D(1)}$, and $e_{k(D)} + d_{k(D)} = o_D(m_{k(D)})$ as $D \to \infty$. Now, set $\sigma := \sigma_*/2$ and define $m(D) := m_{k(D)}$, $e(D) := 4e_{k(D)}$, and $d(D) := d_{k(D)}$. This establishes item $1$ of Assumption~\ref{assumption:mesoscopic}. In the remainder of the proof, take $D$ sufficiently large.

Now, let $\mathcal{Q}_D := (e_{k(D)} \Z) \cap [m(D),D/100]$ and $\mathcal{Z}_D := (e_{k(D)} \Z^2) \cap B_{2D}(0)$. These are $e(D)$-nets of $[m(D),D/100]$ and $B_{2D}(0)$ respectively, which establishes the net requirement in item $3$ of Assumption~\ref{assumption:mesoscopic}.

Finally, $\rho R_{b_{k(D)}} \geq 5D$ implies $B_{4D}(0) \subset B_{\rho R_{b_{k(D)}}}(0)$, so every cell intersecting $B_{4D}(0)$ has diameter at most $d(D)$, establishing item $2$. If $z \in \mathcal{Z}_D$ and $r \in \mathcal{Q}_D$, then $B_{2r}(z) \subset B_{(2 + 1/50)D}(0) \subset B_{\rho R_{b_{k(D)}}}(0)$. Along with the earlier quantum cone scaling argument, this establishes item $3$ of Assumption~\ref{assumption:mesoscopic} for some deterministic $p > 0$.

For item $4$, let $A \ni o$ be finite and connected with $D := \diam(U_A^{\cH})$ sufficiently large. Since $0 \in H_o$, we have $U_A^{\cH} \subset \overline{B_D(0)}$. Proposition~\ref{prop:mated-euclidean-cell-count} with $\delta = 1$ and $k = k(D)$ gives $|A| \leq |\cV_{\mathcal{G}^1}(B_{R_{b_{k(D)}}}(0))| \leq 2^{2k(D)} \leq 4D^{2 \lambda}$. Local finiteness gives $D \to \infty$ uniformly as $|A| \to \infty$, so item $4$ holds with $\nu = 1/(2 \lambda)$ and $c_{\vol} = 4^{-\nu}$.
\end{proof}

Finally, Proposition~\ref{prop:mated-mesoscopic} and Theorem~\ref{thm:poly-sublinear} imply that DLA grown from the root $o$ on the $\gamma$-mated CRT has a.s.\ polynomially sublinear growth. Inspection of the exponent derived in Proposition~\ref{prop:mated-mesoscopic} reveals an inequality which relates the DLA growth bound to the LQG parameter $\gamma$.

\section{DLA on an ergodic scale-free random environment}\label{sec:ergodic}
In this brief section, we sketch the proof of Theorem~\ref{thm:dla-ergodic-env}. Let $\mathcal{M}$ be an embedded lattice as in Definition~\ref{def:ergodic}, and let $G = \mathcal{M}^*$ be its adjacency graph, with root $o = H_0$ and DLA initial cluster $K_0 = \{o\}$. We again assume that harmonic measure from infinity is unique. In this section, good scales are defined using the conductance weighted walk in the obvious way.

\begin{proof}[Proof sketch]
Let $X$ be the conductance weighted walk on $G = \mathcal{M}^*$, with vertex embeddings $\Phi(H) = p_H$ as in Proposition~\ref{prop:ergodic}. Rewritten in the notation of Assumption~\ref{assumption:rw-conv}, \cite[Theorem~3.10]{gwynne2021invarianceprincipleergodicscalefree} gives a.s.\ over $\mathcal{M}$, for every $0 < C < L < \infty$ the convergence
$$
  \sup_{x:\,|\Phi(x)| \leq CR}
  d_{\textrm{Pr}}(
    \cL_x(
      R^{-1} \widehat{Y}^{x;R,L} \mid \mathcal{M}
    ),
    \cL(B_\Sigma^{\Phi(x)/R;L}))
  \to 0
  \text{ as } R \to \infty
$$
where $\widehat{Y}^{x;R,L}$ is the piecewise-linear interpolation through chosen points $p_H$ of the faces, and $B_\Sigma^{z;L}$ is Brownian motion with positive-definite covariance matrix $\Sigma$ started from $z$ and stopped upon leaving $B_L(0)$.

Draw the dual graph as a proper planar map by placing each dual vertex at its chosen face representative and drawing each dual edge inside the union of its two incident faces. By \cite[Lemma 2.9]{gwynne2021invarianceprincipleergodicscalefree}, for every fixed $L < \infty$ we have
\[
  \frac{1}{R}
  \max \left \{
    \diam(H):H \cap B_{LR}(0) \neq \emptyset
  \right \}
  \to 0.
\]
Now, fix $0 < C < L < \infty$ and choose $L' > L$. Set $\delta_R := R^{-1} \max \{ \diam(H): H \cap B_{(L' + 1)R}(0) \neq \emptyset \}$. The preceding limit implies that $\delta_R \to 0$ almost surely. Going forward, we fix an environment on which these almost sure conclusions hold; all walk laws below are conditional on $\mathcal{M}$.

We couple the proper and piecewise-linear embeddings through the same random walk on the graph, started from $x$ with $\Phi(x) \in B_{CR}(0)$, and denote them by $Y$ and $\widehat{Y}$ respectively. For adjacent faces $H,H'$, since $H \cap H' \neq \emptyset$, one has $\diam(H \cup H') \leq \diam(H) + \diam(H')$. It follows that $d_F(R^{-1}Y,R^{-1} \widehat{Y}) \leq 2 \delta_R$, where both paths are stopped when the piecewise-linear interpolation leaves $B_{L'R}(0)$.

Write $B := B_\Sigma^{\Phi(x)/R;L'}$. Since the Prokhorov metric is induced by $1 \wedge d_F$, the above implies
$$d_{\textrm{Pr}}(\cL_x(R^{-1}Y),
\cL_x(R^{-1} \widehat{Y})) \leq 2 \delta_R.$$ By the triangle inequality,
$d_{\textrm{Pr}}(\cL_x(R^{-1}Y), \cL(B)) \leq 2 \delta_R + d_{\textrm{Pr}}(\cL_x(R^{-1} \widehat{Y}), \cL(B))$. The second term tends to zero uniformly over starting vertices in $B_{CR}(0)$ by Proposition~\ref{prop:ergodic}, while $\delta_R \to 0$ is already shown.

Finally, we deal with a slight mismatch in stopping times between the two walks, since the above $Y$ is stopped when $\widehat{Y}$ exits $B_{L'R}(0)$; we prove that $Y$ stopped when it exits $B_{LR}(0)$, after multiplication by $R^{-1}$, converges to Brownian motion stopped when exiting $B_L(0)$. The terminal point of $Y$ has modulus at least $(L' - 2 \delta_R)R > LR$ for all large $R$, so this exit occurs before the common graph-time interval ends.

Fix a sequence $R_n \to \infty$ and starting vertices $x_n$ with $\Phi(x_n) \in B_{CR_n}(0)$, and suppose first that $d_{\textrm{Pr}}(\cL_{x_n}(R_n^{-1}Y), \cL(B_\Sigma^{z;L'})) \to 0$ for some $z \in \overline{B_C(0)}$. By the Skorokhod representation theorem, we may realize this convergence as a.s.\ convergence under $d_F$; denote by $\gamma_n$ the corresponding curves, with their parametrizations chosen so that $\gamma_n \to \gamma$ uniformly on $[0,1]$. Here $\gamma$ is the curve $B_\Sigma^{z;L'}$ parametrized by $[0,1]$ in the obvious way (so $\gamma$ exits $B_{L'}(0)$ at $t = 1$). Let $\tau_n$ and $\tau$ be their respective first entry into $\C \setminus B_L(0)$. For each sufficiently small $\delta > 0$, the compact set $\gamma([0, \tau - \delta])$ lies a positive distance inside $B_L(0)$, so $\tau_n \geq \tau - \delta$ for all sufficiently large $n$. On the other hand, a.s.\ there exists $s \in (\tau, \tau + \delta)$ with $|\gamma(s)| > L$. Uniform convergence then gives $|\gamma_n(s)| > L$, and so $\tau_n \leq s < \tau + \delta$, for all sufficiently large $n$, which implies $\tau_n \to \tau$ a.s. By continuity of $\gamma$, the parametrizations $u \mapsto \gamma_n(u \tau_n)$ converge uniformly to $u \mapsto \gamma(u \tau)$, so the corresponding initial curve segments converge in $d_F$. Thus the rescaled random walk under $\Phi$, stopped upon leaving $B_{LR_n}(0)$, converges weakly to $B_\Sigma^{z;L}$.

To obtain uniformity over starting vertices, suppose that uniform convergence fails along some sequence $R_n \to \infty$ and $x_n$, with some $C < L$, with the Prokhorov distances bounded away from zero. By sequential compactness, we can choose a convergent subsequence $z_n := \Phi(x_n)/R_n \to z \in \overline{B_C(0)}$. Write the Brownian motions as $z_n + W$ and $z + W$, where $W$ is the same Brownian motion of covariance $\Sigma$ started at $0$. Their distance before stopping is $|z_n - z| \to 0$. The same $\tau - \delta, \tau + \delta$ comparison above, applied at each radius $a \in \{L,L'\}$, shows that their departure times and stopped curves converge, so $d_{\textrm{Pr}}(\cL(B_\Sigma^{z_n;a}), \cL(B_\Sigma^{z;a})) \to 0$. Now, the (non-uniform) convergence $d_{\textrm{Pr}}(\cL_{x_n}(R_n^{-1}Y \mid \mathcal{M}), \cL(B_\Sigma^{\Phi(x_n)/R_n;L'})) \to 0$ and $d_{\textrm{Pr}}(\cL(B_\Sigma^{z_n;L'}), \cL(B_\Sigma^{z;L'})) \to 0$, along with the triangle inequality show that $d_{\textrm{Pr}}(\cL_{x_n}(R_n^{-1}Y \mid \mathcal{M}), \cL(B_\Sigma^{z;L'})) \to 0$. The preceding paragraph, together with $a = L$ and the triangle inequality, show that the Prokhorov distances tend to zero, which is a contradiction.

Lastly, we deal with the anisotropy of $\Sigma$. Let $\eta$ be as in Claim~\ref{claim:brownian-robust-intersection} and let $U_\theta$ denote a rotation by $\theta$ and set $\Sigma_\theta := U_{-\theta} \Sigma U_\theta$. Write $\P_z^\Sigma$ for the law of Brownian motion $\widetilde B$ started at $z$ with covariance $\Sigma$. If $\Gamma_\rho$ is the polygonal path in Figure~\ref{fig:polygonal-path} with initial radius $\rho \in [1 - \eta,1 + \eta]$, set $\eta_0 := \eta/2$ and define
$$
  p := \frac{1}{2}
  \inf_{\substack{\theta \in [0,2 \pi]\\
                  \rho \in [1 - \eta,1 + \eta]}}
  \P_{(\rho,0)}^{\Sigma_\theta}
  [\exists \, t > 0:\ d_F(\widetilde B|_{[0,t]}, \Gamma_\rho) < \eta_0]
$$
which is strictly positive by compactness and continuity.

The good scales argument now applies identically. For every $K$, there exists an $N(K) < \infty$ such that for every finite, connected, rooted set $A$ with $|A| > N(K)$ and every $v \in A$, there are $(K,p)$-good scales for $(A,v)$. This implies
\[
  \harm_A^\infty(\mathcal{E}_A(v))
  \leq (1 - p)^K.
\]
The argument for DLA growth then follows identically to Proposition~\ref{prop:helper-diam}, and a.s.\ sublinearity of the internal diameter follows for DLA grown from $H_0$.
\end{proof}

\section{Open problems}\label{sec:open-probs}
We briefly mention two natural directions of further research. The first regards the identification of the growth exponent of the internal diameter of DLA. Theorem~\ref{thm:poly-sublinear} bounds the growth exponent of the internal diameter of DLA. As mentioned in Section~\ref{sec:intro}, the paper \cite{gwynne2021external} provides a result for the growth exponent of the external diameter of DLA; this result is rather strong, as it identifies precisely the growth exponent, but it is limited by the special qualities of the UITM. A natural next question is the following: 
\begin{problem}
    With $(G,o)$ a whole-plane $\gamma$-mated CRT rooted at $o = 0$, and $(K_n)_{n \in \N_0}$ a DLA process from $o$ to $\infty$, show that
    $$
    \lim_{n \to \infty} \frac{\log \diam(K_n)}{\log n}
    $$
    exists, and identify its value.
\end{problem}
The above problem is essentially \cite[Problem 8.1]{gwynne2021external} stated for the $\gamma$-mated CRT.

The next more open ended problem deals with a limitation of the techniques of this paper. To deduce a.s.\ sublinear growth of the DLA process on an RPM, we require Assumption~\ref{assumption:rw-conv} or \ref{assumption:mesoscopic} to hold a.s.\ for the relevant model. As demonstrated in this paper, this does indeed hold for an extremely general class of random planar geometries. However, it does fail to capture the behavior of some natural models which do exhibit slightly weaker invariance properties; in particular, the recent work \cite{holden2026tutte} exhibits quenched convergence \textit{in probability} of the random walk (uniform on distinct neighboring vertices) on finite spanning-tree-weighted maps to Brownian motion, and indeed the work of \cite{gwynne2021external} shows that the UITM is an especially natural geometry on which to study DLA. As we will discuss, these in-probability invariance theorems alone are insufficient to prove sublinearity of the internal diameter of DLA. Thus, the following problem is a natural modification, and asks if further assumptions are sufficient to recover a.s.\ sublinear diameter growth.

\begin{problem}
    Adapt the methods of the current paper to RPMs satisfying an in probability quenched invariance principle under further assumptions.
\end{problem}

To observe that in-probability convergence alone cannot yield the sublinear bound on DLA, consider the following example. Set $r_k = 4^{2^k}$ and $Q_r = [-r,r]^2 \cap \Z^2$ and independently with probability $1/k$ replace $Q_{2r_k} \setminus Q_{r_k}$ in $\Z^2$ by a single nearest-neighbor path with $L_k = r_k(2r_k + 1)$ vertices, retaining only consecutive edges and the two connections to the inner and outer regions. The resulting graph is connected, recurrent, and properly embedded in the plane, with maximum degree four and unique harmonic measure from infinity. Set $o = 0$ and take $K_0 = \{0\}$. For each fixed $L > 0$, with probability tending to one all modifications relevant to a walk stopped on exiting $B_{LR}(0)$ lie inside $B_{4R^{1/4}}(0)$. The lattice invariance principle and planar annulus hitting estimate \cite[Theorem 3.4.2 and Proposition 6.4.1]{lawler2010random} give the convergence in Assumption~\ref{assumption:rw-conv} in probability. However, infinitely many modifications occur a.s.\ and DLA must traverse them consecutively. At their completion times $N_k$, we have $N_k \leq (2r_k + 1)^2 + L_k$ and $\diam(K_{N_k}) \geq L_k$, so
$$
  \limsup_{n \to \infty} \frac{\diam(K_n)}{n} \geq \frac{1}{3}
$$
a.s.\ For $a_m(G,o)$, monotonicity makes convergence to zero in probability equivalent to a.s.\ convergence, so Proposition~\ref{prop:helper-diam} gives sublinear DLA diameter jointly a.s.\ However, in the counterexample we have $a_m(G,0) = 1$ for every $m$ a.s., since arbitrarily large rooted sets can end inside a long path with a single entrance from infinity. Likewise, convergence in probability of the walk-law discrepancy uniformly over all later scales recovers a.s.\ convergence.

Specifically, for the UITM it seems likely that one may be able to achieve analogous statements to \cite{gwynne2019harmonicfunctionsmatedcrtmaps} and therefore achieve the a.s. Assumption~\ref{assumption:mesoscopic} due to the recent work \cite{holden2026tutte}.

\begin{problem}
    Using \cite{holden2026tutte}, establish analogues of the quantitative quenched annular estimates of \cite{gwynne2019harmonicfunctionsmatedcrtmaps} on the UITM, and then verify Assumption~\ref{assumption:mesoscopic} for the UITM.
\end{problem}

\bibliographystyle{plain}
\bibliography{refs}

@article{AB99,
    title = "H{\"o}lder regularity and dimension bounds for random curves",
    author = "M. Aizenman and A. Burchard",
    year = "1999",
    doi = "10.1215/S0012-7094-99-09914-3",
    language = "English (US)",
    volume = "99",
    pages = "419--453",
    journal = "Duke Mathematical Journal",
    issn = "0012-7094",
    publisher = "Duke University Press",
    number = "3",
}

@article{amir2016one,
  title={One-dimensional long-range diffusion-limited aggregation I},
  author={Amir, Gideon and Angel, Omer and Benjamini, Itai and Kozma, Gady},
  journal={The Annals of Probability},
  pages={3546--3579},
  year={2016},
  publisher={JSTOR}
}

@article{ang2026quantum,
  title={Quantum Loewner evolution in quantum natural time: phases and Markov properties},
  author={Ang, Morris and Mithal, Deven},
  journal={arXiv preprint arXiv:2605.03385},
  year={2026}
}

@article{asselah2013sublogarithmic,
  author  = {Asselah, Amine and Gaudilli\`ere, Alexandre},
  title   = {Sublogarithmic fluctuations for internal {DLA}},
  journal = {The Annals of Probability},
  year    = {2013},
  volume  = {41},
  number  = {3A},
  pages   = {1160--1179},
  doi     = {10.1214/11-AOP735}
}

@article{benjamini2008diffusion,
  title={Diffusion limited aggregation on a cylinder},
  author={Benjamini, Itai and Yadin, Ariel},
  journal={Communications in Mathematical Physics},
  volume={279},
  number={1},
  pages={187--223},
  year={2008},
  publisher={Springer}
}

@article{benjamini2017upper,
  title={Upper bounds on the growth rate of diffusion limited aggregation},
  author={Benjamini, Itai and Yadin, Ariel},
  journal={arXiv preprint arXiv:1705.06095},
  year={2017}
}

@article{BenjaminiYadinSpectralGap,
  author  = {Benjamini, Itai and Yadin, Ariel},
  title   = {Harmonic measure in the presence of a spectral gap},
  journal = {Annales de l'Institut Henri Poincar\'e, Probabilit\'es et Statistiques},
  year    = {2016},
  volume  = {52},
  number  = {3},
  pages   = {1050--1060},
  doi     = {10.1214/15-AIHP670}
}

@article{berestycki2024harnack,
  title={Harnack inequality and one-endedness of UST on reversible random graphs: N. Berestycki, D. van Engelenburg},
  author={Berestycki, Nathana{\"e}l and van Engelenburg, Diederik},
  journal={Probability Theory and Related Fields},
  volume={188},
  number={1},
  pages={487--548},
  year={2024},
  publisher={Springer}
}

@article{berger2022growth,
  title={Growth of stationary hastings--levitov},
  author={Berger, Noam and Procaccia, Eviatar B and Turner, Amanda},
  journal={The Annals of Applied Probability},
  volume={32},
  number={5},
  pages={3331--3360},
  year={2022},
  publisher={Institute of Mathematical Statistics}
}

@article{bubeck2015asymptotic,
  title={Asymptotic behavior of the Eden model with positively homogeneous edge weights},
  author={Bubeck, S{\'e}bastien and Gwynne, Ewain},
  journal={arXiv preprint arXiv:1508.05140},
  year={2015}
}

@misc{bonannoFormalizations,
  author = {Bonanno, Leonardo and Mithal, Deven},
  title = {DLA-internal-Diameter},
  howpublished = {GitHub repository, \url{https://github.com/liouvillequantumguy/Formalizations/tree/main/DLA-internal-diameter}}
}

@article{bou2022dynamic,
  title={Dynamic dimensional reduction in the Abelian sandpile},
  author={Bou-Rabee, Ahmed},
  journal={Communications in Mathematical Physics},
  volume={390},
  number={2},
  pages={933--958},
  year={2022},
  publisher={Springer}
}

@article{bou2024shape,
  title={A shape theorem for exploding sandpiles},
  author={Bou-Rabee, Ahmed},
  journal={The Annals of Applied Probability},
  volume={34},
  number={1A},
  pages={714--742},
  year={2024},
  publisher={Institute of Mathematical Statistics}
}

@article{bou2024integer,
  title={Integer superharmonic matrices on the F-lattice},
  author={Bou-Rabee, Ahmed},
  journal={Advances in Mathematics},
  volume={436},
  pages={109400},
  year={2024},
  publisher={Elsevier}
}

@article{bou2024internal,
  title={Internal DLA on mated-CRT maps},
  author={Bou-Rabee, Ahmed and Gwynne, Ewain},
  journal={The Annals of Probability},
  volume={52},
  number={6},
  pages={2173--2237},
  year={2024},
  publisher={Institute of Mathematical Statistics}
}

@article{bou2021convergence,
  title={Convergence of the random Abelian sandpile},
  author={Bou-Rabee, Ahmed},
  journal={The Annals of Probability},
  volume={49},
  number={6},
  pages={3168--3196},
  year={2021},
  publisher={JSTOR}
}

@article{cox1981some,
  title={Some limit theorems for percolation processes with necessary and sufficient conditions},
  author={Cox, J Theodore and Durrett, Richard},
  journal={The Annals of Probability},
  pages={583--603},
  year={1981},
  publisher={JSTOR}
}

@article{curien2017scaling,
  author  = {Curien, Nicolas and Le Gall, Jean-Fran\c{c}ois},
  title   = {Scaling limits for the peeling process on random maps},
  journal = {Annales de l'Institut Henri Poincar{\'e}, Probabilit{\'e}s et Statistiques},
  volume  = {53},
  number  = {1},
  pages   = {322--357},
  year    = {2017},
  doi     = {10.1214/15-AIHP718}
}

@inproceedings{curien2019first,
  title={First-passage percolation and local modifications of distances in random triangulations},
  author={Curien, Nicolas and Le Gall, Jean-Fran{\c{c}}ois},
  booktitle={Annales scientifiques de l'{\'E}cole Normale Sup{\'e}rieure},
  volume={52},
  number={3},
  pages={631--701},
  year={2019}
}

@article{carleson2001aggregation,
  title={Aggregation in the plane and Loewner's equation},
  author={Carleson, Lennart and Makarov, Nikolai},
  journal={Communications in Mathematical Physics},
  volume={216},
  number={3},
  pages={583--607},
  year={2001},
  publisher={Springer}
}

@article{duplantier2011liouville,
  title={Liouville quantum gravity and KPZ},
  author={Duplantier, Bertrand and Sheffield, Scott},
  journal={Inventiones mathematicae},
  volume={185},
  number={2},
  pages={333--393},
  year={2011},
  publisher={Springer}
}

@article{duplantier2021liouville,
  title={Liouville quantum gravity as a mating of trees},
  author={Duplantier, Bertrand and Miller, Jason R and Sheffield, Scott},
  journal={Ast{\'e}risque},
  volume={427},
  year={2021}
}

@inproceedings{eden1961two,
  title={A two-dimensional growth process},
  author={Eden, Murray and others},
  booktitle={Proceedings of the fourth Berkeley symposium on mathematical statistics and probability},
  volume={4},
  pages={223--239},
  year={1961},
  organization={Biology and Problems of Health, University of California Press, Berkeley, CA}
}

@article{gwynne2019distance,
  title={A distance exponent for Liouville quantum gravity},
  author={Gwynne, Ewain and Holden, Nina and Sun, Xin},
  journal={Probability Theory and Related Fields},
  volume={173},
  number={3},
  pages={931--997},
  year={2019},
  publisher={Springer}
}

@article{gwynne2019harmonicfunctionsmatedcrtmaps,
  author  = {Gwynne, Ewain and Miller, Jason and Sheffield, Scott},
  title   = {Harmonic functions on mated-{CRT} maps},
  journal = {Electronic Journal of Probability},
  year    = {2019},
  volume  = {24},
  pages   = {Paper No. 58, 55},
  doi     = {10.1214/19-EJP325}
}

@article{gwynne2020anomalous,
  title={Anomalous diffusion of random walk on random planar maps},
  author={Gwynne, Ewain and Hutchcroft, Tom},
  journal={Probability Theory and Related Fields},
  volume={178},
  number={1},
  pages={567--611},
  year={2020},
  publisher={Springer}
}

@article{gwynne2021external,
  author  = {Gwynne, Ewain and Pfeffer, Joshua},
  title   = {External diffusion-limited aggregation on a spanning-tree-weighted random planar map},
  journal = {The Annals of Probability},
  year    = {2021},
  volume  = {49},
  number  = {4},
  pages   = {1633--1676},
  doi     = {10.1214/20-AOP1486}
}

@article{gwynne2021invarianceprincipleergodicscalefree,
  author  = {Gwynne, Ewain and Miller, Jason and Sheffield, Scott},
  title   = {An invariance principle for ergodic scale-free random environments},
  journal = {Acta Mathematica},
  year    = {2022},
  volume  = {228},
  number  = {2},
  pages   = {303--384},
  doi     = {10.4310/ACTA.2022.v228.n2.a2}
}

@article{gwynne2021tutteembeddingmatedcrtmap,
  author  = {Gwynne, Ewain and Miller, Jason and Sheffield, Scott},
  title   = {The {Tutte} embedding of the mated-{CRT} map converges to {Liouville} quantum gravity},
  journal = {The Annals of Probability},
  year    = {2021},
  volume  = {49},
  number  = {4},
  pages   = {1677--1717},
  doi     = {10.1214/20-AOP1487}
}

@article{hastings2001fractal,
  title={Fractal to nonfractal phase transition in the dielectric breakdown model},
  author={Hastings, MB},
  journal={Physical Review Letters},
  volume={87},
  number={17},
  pages={175502},
  year={2001},
  publisher={APS}
}

@article{hastings1998laplacian,
  title={Laplacian growth as one-dimensional turbulence},
  author={Hastings, Matthew B and Levitov, Leonid S},
  journal={Physica D: Nonlinear Phenomena},
  volume={116},
  number={1-2},
  pages={244--252},
  year={1998},
  publisher={Elsevier}
}

@article{HiggsALE,
  author  = {Higgs, Frankie},
  title   = {{SLE} scaling limits for a {Laplacian} random growth model},
  journal = {Annales de l'Institut Henri Poincar\'e, Probabilit\'es et Statistiques},
  year    = {2022},
  volume  = {58},
  number  = {3},
  pages   = {1712--1739},
  doi     = {10.1214/21-AIHP1217}
}

@article{holden2026tutte,
  title={The Tutte embedding of tree-weighted planar maps converges to Liouville quantum gravity},
  author={Holden, Nina},
  journal={arXiv preprint arXiv:2610.10501},
  year={2026}
}

@inproceedings{hutchcroft2019uniform,
  title={Uniform spanning forests of planar graphs},
  author={Hutchcroft, Tom and Nachmias, Asaf},
  booktitle={Forum of Mathematics, Sigma},
  volume={7},
  pages={e29},
  year={2019},
  organization={Cambridge University Press}
}

@inproceedings{hammersley1965first,
  title={First-passage percolation, subadditive processes, stochastic networks, and generalized renewal theory},
  author={Hammersley, John M and Welsh, Dominic JA},
  booktitle={Bernoulli 1713, Bayes 1763, Laplace 1813: Anniversary Volume. Proceedings of an International Research Seminar Statistical Laboratory University of California, Berkeley 1963},
  pages={61--110},
  year={1965},
  organization={Springer}
}

@article{jerison2012logarithmic,
  title={Logarithmic fluctuations for internal DLA},
  author={Jerison, David and Levine, Lionel and Sheffield, Scott},
  journal={Journal of the American Mathematical Society},
  volume={25},
  number={1},
  pages={271--301},
  year={2012}
}

@article{kesten1987long,
  title={How long are the arms in DLA?},
  author={Kesten, Harry},
  journal={Journal of Physics A: Mathematical and General},
  volume={20},
  number={1},
  pages={L29--L33},
  year={1987}
}

@article{kesten1987hitting,
  title={Hitting probabilities of random walks on Zd},
  author={Kesten, Harry},
  journal={Stochastic Processes and their Applications},
  volume={25},
  pages={165--184},
  year={1987},
  publisher={Elsevier}
}

@article{kesten1990upper,
  title={Upper bounds for the growth rate of DLA},
  author={Kesten, Harry},
  journal={Physica A: Statistical Mechanics and its Applications},
  volume={168},
  number={1},
  pages={529--535},
  year={1990},
  publisher={Elsevier}
}

@book{lawler2008conformally,
       AUTHOR = {Lawler, Gregory F.},
     TITLE = {Conformally invariant processes in the plane},
    SERIES = {Mathematical Surveys and Monographs},
    VOLUME = {114},
 PUBLISHER = {American Mathematical Society},
   ADDRESS = {Providence, RI},
      YEAR = {2005},
     PAGES = {xii+242},
      ISBN = {0-8218-3677-3},
   MRCLASS = {60-02 (30-02 30C35 31A15 60J65 81T40 82B27)},
  MRNUMBER = {2129588 (2006i:60003)},
MRREVIEWER = {Zhen-Qing Chen},
}

@book{lawler2010random,
  author    = {Lawler, Gregory F. and Limic, Vlada},
  title     = {Random Walk: A Modern Introduction},
  series    = {Cambridge Studies in Advanced Mathematics},
  volume    = {123},
  publisher = {Cambridge University Press},
  address   = {Cambridge},
  year      = {2010},
  doi       = {10.1017/CBO9780511750854}
}

@article{lawler1992internal,
  title={Internal diffusion limited aggregation},
  author={Lawler, Gregory F and Bramson, Maury and Griffeath, David},
  journal={The Annals of Probability},
  pages={2117--2140},
  year={1992},
  publisher={JSTOR}
}

@article{lawler2004beurling,
  author  = {Lawler, Gregory F. and Limic, Vlada},
  title   = {The {Beurling} estimate for a class of random walks},
  journal = {Electronic Journal of Probability},
  year    = {2004},
  volume  = {9},
  number  = {27},
  pages   = {846--861},
  doi     = {10.1214/EJP.v9-228}
}

@article{levine2010scaling,
  title={Scaling limits for internal aggregation models with multiple sources},
  author={Levine, Lionel and Peres, Yuval},
  journal={Journal d'Analyse Math{\'e}matique},
  volume={111},
  number={1},
  pages={151--219},
  year={2010},
  publisher={Springer}
}

@article{losev2025long,
  title={How Long are the Arms in DBM?},
  author={Losev, Ilya and Smirnov, Stanislav},
  journal={Communications in Mathematical Physics},
  volume={406},
  number={4},
  pages={93},
  year={2025},
  publisher={Springer}
}

@article{lawler2011conformal,
    AUTHOR = {Lawler, Gregory F. and Schramm, Oded and Werner, Wendelin},
     TITLE = {Conformal invariance of planar loop-erased random walks and uniform spanning trees},
   JOURNAL = {Ann. Probab.},
  FJOURNAL = {The Annals of Probability},
    VOLUME = {32},
      YEAR = {2004},
    NUMBER = {1B},
     PAGES = {939--995},
      ISSN = {0091-1798},
     CODEN = {APBYAE},
   MRCLASS = {82B41 (60G50)},
  MRNUMBER = {2044671 (2005f:82043)},
MRREVIEWER = {Olivier Raimond},
       DOI = {10.1214/aop/1079021469},
       URL = {http://dx.doi.org/10.1214/aop/1079021469},
       eprint={\arxiv{math/0112234}}
}

@misc{ms-equivalence,
      title={{Liouville quantum gravity and the Brownian map I: The QLE(8/3,0) metric}}, 
      author={Jason Miller and Scott Sheffield},
      year={2019},
      eprint={1507.00719},
      archivePrefix={arXiv},
      primaryClass={math.PR},
      url={https://arxiv.org/abs/1507.00719}, 
}

@article{ms-equiavlence-2,
    AUTHOR = {Miller, Jason and Sheffield, Scott},
     TITLE = {Liouville quantum gravity and the {B}rownian map {II}:
              {G}eodesics and continuity of the embedding},
   JOURNAL = {Ann. Probab.},
  FJOURNAL = {The Annals of Probability},
    VOLUME = {49},
      YEAR = {2021},
    NUMBER = {6},
     PAGES = {2732--2829},
      ISSN = {0091-1798},
   MRCLASS = {60D05 (60J67 60J80 83C45)},
  MRNUMBER = {4348679},
       DOI = {10.1214/21-aop1506},
       URL = {https://doi.org/10.1214/21-aop1506},
   eprint = {\arxiv{1605.03563}}, 
}

@article{MJ02,
  author  = {Mathiesen, Joachim and Jensen, Mogens H.},
  title   = {Tip splittings and phase transitions in the dielectric breakdown model: Mapping to the diffusion-limited aggregation model},
  journal = {Physical Review Letters},
  volume  = {88},
  number  = {23},
  pages   = {235505},
  year    = {2002},
  doi     = {10.1103/PhysRevLett.88.235505}
}

@article{MJB08,
  author  = {Mathiesen, Joachim and Jensen, Mogens H. and Bakke, Jan {\O{}}ystein Haavig},
  title   = {Dimensions, maximal growth sites, and optimization in the dielectric breakdown model},
  journal = {Physical Review E},
  volume  = {77},
  number  = {6},
  pages   = {066203},
  year    = {2008},
  doi     = {10.1103/PhysRevE.77.066203}
}

@book{morters2010brownian,
  author    = {M{\"o}rters, Peter and Peres, Yuval},
  title     = {Brownian Motion},
  series    = {Cambridge Series in Statistical and Probabilistic Mathematics},
  volume    = {30},
  publisher = {Cambridge University Press},
  address   = {Cambridge},
  year      = {2010},
  doi       = {10.1017/CBO9780511750489}
}

@article{ms-qle,
   title={{Quantum Loewner evolution}},
   volume={165},
   ISSN={0012-7094},
   url={http://dx.doi.org/10.1215/00127094-3627096},
   DOI={10.1215/00127094-3627096},
   number={17},
   journal={Duke Mathematical Journal},
   publisher={Duke University Press},
   author={Miller, Jason and Sheffield, Scott},
   year={2016},
   pages = {3241--3378},
   month=nov 
}

@article{NorrisSilvestriTurnerALE,
  author  = {Norris, James and Silvestri, Vittoria and Turner, Amanda},
  title   = {Stability of regularized {Hastings--Levitov} aggregation in the subcritical regime},
  journal = {Communications in Mathematical Physics},
  year    = {2024},
  volume  = {405},
  number  = {3},
  pages   = {74},
  doi     = {10.1007/s00220-024-04960-5}
}

@article{NorrisTurnerHL,
  author  = {Norris, James and Turner, Amanda},
  title   = {{Hastings--Levitov} aggregation in the small-particle limit},
  journal = {Communications in Mathematical Physics},
  year    = {2012},
  volume  = {316},
  number  = {3},
  pages   = {809--841},
  doi     = {10.1007/s00220-012-1552-6}
}

@article{niemeyer1984dbm,
  title={Fractal dimension of dielectric breakdown},
  author={Niemeyer, Lucian and Pietronero, Luciano and Wiesmann, Hans J},
  journal={Physical Review Letters},
  volume={52},
  number={12},
  pages={1033--1036},
  year={1984},
  publisher={APS}
}

@article{ProcacciaYeZhangStationary,
  author  = {Procaccia, Eviatar B. and Ye, Jiayan and Zhang, Yuan},
  title   = {Stationary {DLA} is well defined},
  journal = {Journal of Statistical Physics},
  year    = {2020},
  volume  = {181},
  number  = {4},
  pages   = {1089--1111},
  doi     = {10.1007/s10955-020-02619-8}
}

@article{ProcacciaZhangStationary,
  author  = {Procaccia, Eviatar B. and Zhang, Yuan},
  title   = {Stationary harmonic measure and {DLA} in the upper half plane},
  journal = {Journal of Statistical Physics},
  year    = {2019},
  volume  = {176},
  number  = {4},
  pages   = {946--980},
  doi     = {10.1007/s10955-019-02327-y}
}

@article{Ric73,
  author  = {Richardson, Daniel},
  title   = {Random growth in a tessellation},
  journal = {Proceedings of the Cambridge Philosophical Society},
  volume  = {74},
  number  = {3},
  pages   = {515--528},
  year    = {1973},
  doi     = {10.1017/S0305004100077288}
}

@article{RZ05,
  author  = {Rohde, Steffen and Zinsmeister, Michel},
  title   = {Some remarks on {Laplacian} growth},
  journal = {Topology and its Applications},
  volume  = {152},
  number  = {1--2},
  pages   = {26--43},
  year    = {2005},
  doi     = {10.1016/j.topol.2004.08.013}
}

@article{schramm0,
    AUTHOR = {Schramm, Oded},
     TITLE = {Scaling limits of loop-erased random walks and uniform spanning trees},
   JOURNAL = {Israel J. Math.},
  FJOURNAL = {Israel Journal of Mathematics},
    VOLUME = {118},
      YEAR = {2000},
     PAGES = {221--288},
      ISSN = {0021-2172},
     CODEN = {ISJMAP},
   MRCLASS = {60K35 (30C85 60H15 82B27 82B44)},
  MRNUMBER = {1776084 (2001m:60227)},
MRREVIEWER = {Almut Burchard},
       DOI = {10.1007/BF02803524},
       URL = {http://dx.doi.org/10.1007/BF02803524},
       eprint={\arxiv{math/9904022}}
}

@article{smirnov-cardy,
AUTHOR = {Smirnov, Stanislav},
     TITLE = {Critical percolation in the plane: conformal invariance, {C}ardy's formula, scaling limits},
   JOURNAL = {C. R. Acad. Sci. Paris S\'er. I Math.},
  FJOURNAL = {Comptes Rendus de l'Acad\'emie des Sciences. S\'erie I. Math\'ematique},
    VOLUME = {333},
      YEAR = {2001},
    NUMBER = {3},
     PAGES = {239--244},
      ISSN = {0764-4442},
     CODEN = {CASMEI},
   MRCLASS = {60K35 (82B43)},
  MRNUMBER = {1851632 (2002f:60193)},
       DOI = {10.1016/S0764-4442(01)01991-7},
       URL = {http://dx.doi.org/10.1016/S0764-4442(01)01991-7},
eprint={\arxiv{0909.4499}}
}

@article{smirnov2010conformal,
  title={Conformal invariance in random cluster models. I. Holmorphic fermions in the Ising model},
  author={Smirnov, Stanislav},
  journal={Annals of mathematics},
  pages={1435--1467},
  year={2010},
  publisher={JSTOR}
}

@article{SolaTurnerViklund,
  author  = {Sola, Alan and Turner, Amanda and Viklund, Fredrik},
  title   = {One-dimensional scaling limits in a planar {Laplacian} random growth model},
  journal = {Communications in Mathematical Physics},
  year    = {2019},
  volume  = {371},
  number  = {1},
  pages   = {285--329},
  doi     = {10.1007/s00220-019-03460-1}
}

@article{witten1981diffusion,
  title={Diffusion-limited aggregation, a kinetic critical phenomenon},
  author={Witten Jr, Thomas A and Sander, Leonard M},
  journal={Physical review letters},
  volume={47},
  number={19},
  pages={1400},
  year={1981},
  publisher={APS}
}

@article{witten1983diffusion,
  author  = {Witten, Thomas A. and Sander, Leonard M.},
  title   = {Diffusion-limited aggregation},
  journal = {Physical Review B},
  year    = {1983},
  volume  = {27},
  number  = {9},
  pages   = {5686--5697},
  doi     = {10.1103/PhysRevB.27.5686}
}

\end{document}